\documentclass[a4paper,12pt]{article}
\usepackage[dvips]{epsfig}
\usepackage{amsmath,amssymb,euscript,graphicx,amsfonts,pifont}
\usepackage{enumerate,color}
\usepackage{graphicx}
\usepackage{hyperref}
\usepackage{tikz}
\usetikzlibrary{calc, decorations.markings, arrows}

\newtheorem{theorem}{Theorem}[section]
\newtheorem{proposition}[theorem]{Proposition}
\newtheorem{lem}[theorem]{Lemma}

\newcounter{case}

\renewcommand{\thecase}{\arabic{case}}

\newcounter{subcase}

\numberwithin{subcase}{case}

\def\a{\alpha} \def\b{\beta} \def\g{\gamma} \def\d{\delta} \def\e{\varepsilon}
  \def\t{\tau}  
\def\th{\theta} \def\ld{\lambda} \def\ph{\varphi}
\def\Del{\Delta}  
 \def\O{\Omega} \def\G{\Gamma}
\def\D{{\rm \Delta }}

\def\o{\overline}  
\def\di{\bigm|} \def\lg{\langle} \def\rg{\rangle}

\def\PSL{\hbox{\rm PSL}}  
\def\Aut{\hbox{\rm Aut\,}}  
 \def\soc{\hbox{\rm soc}} 
  \def\mod{\hbox{\rm mod }}
  
 \def\PG{\hbox{\rm PG}} \def\PGL{\hbox{\rm PGL}}
  \def\GL{\hbox{\rm GL}}  \def\P\GL{\hbox{\rm P\GL}}
  \def\FF{{\hbox{\sf F\kern-.43emF}}}
\def\cal{\mathcal}  
\def\PGammaL{{\rm P\Gamma L}}  \def\PSp{{\rm PSp}}
\def\Sp{{\rm Sp}} \def\PSU{{\rm PSU}} \def\SU{{\rm SU}}   \def\SO{{\rm SO}}
\def\GO{{\rm GO}} \def\diag{{\rm diag}}

\def\orho{\overline{\rho}}

\def\P{\hbox{\rm P}}

\def\rad{\hbox{\rm rad}}

\def\ZZ{\mathbb{Z}} 

\def\nd{\mathrel{\bigm|\kern-.7em/}} 
 \def\f{\noindent}
\def\qed{\hfill $\Box$} \def\demo{\f {\bf Proof}\hskip10pt}

\begin{document}

	\begin{center}
		{\bf\large Hamilton Cycles of
			Graphs  of Rank 3 \\ Classical Group Cases}
		\footnote{This work was supported in part  by the National Natural Science Foundation of China
			(12471332).}
	\end{center}
	
	
	\begin{center}
		Shaofei Du and Chen Yao\footnote{Corresponding author: 2240501007@cnu.edu.cn}\\
		\medskip
		{\it {\small
				School of Mathematical Sciences, Capital Normal University,\\
				Bejing 100048, People's Republic of China\\
		}}
	\end{center}

	\renewcommand{\thefootnote}{\empty}

\begin{abstract} This paper contributes to a long standing Lov\'asz problem regarding existence of Hamilton paths in vertex-transitive graphs.
 A connected simple graph is said to be of rank $r$ if  its automorphism group is a transitive group of   rank $r$ on vertices.
 Since  every  graph of rank 2  is complete,  a basic step is to  determine  the Hamiltonicity of graphs  of  rank 3.
  It is ease to see that every   graph  of rank 3  whose  automorphism group  is imprimitive must contain a  Hamilton cycle.
 In this paper,   it will be shown that with the exception  of the Petersen graph,  every connected graph of rank 3 contains a  Hamilton cycle,
  provided its automorphism group is  a  primitive group of rank 3,  whose socle is a classical  group.
 The remaining  primitive  groups of rank 3  will be discussed   in our other papers.
 \end{abstract}

\section{Introduction}
\label{sec:intro}
\indent
A simple path (reps. cycle) going through all vertices of a graph is called a Hamilton path (resp. cycle).
In 1969, Lov\'asz \cite{L70} asked that if there exists  a
finite  connected vertex-transitive graph of order more than 3 without a Hamilton path; and   in 1981, Alspach  \cite{A81} asked if there exists an infinite number of connected vertex-transitive graphs that do not have a Hamiltonian cycle.
Till now,  with the exception of $K_2$, only four connected
vertex-transitive graphs that do not have a Hamilton cycle are known to exist.
These four graphs are the Petersen graph, the Coxeter graph and the two graphs obtained
from them by replacing each vertex by a triangle.
The fact that none of these four graphs
is a Cayley graph  has led to a folklore conjecture that except for $K_2$,
every Cayley graph has a Hamilton cycle.

The existence of Hamilton paths, and in some cases
also Hamilton cycles, in  connected vertex-transitive graphs
has been shown for graphs of particular orders, such as, $kp$, where $k\le 6$,
$p^j$, where $j\le 5$, and $2p^2$, where throughout this paper $p$  denotes a prime, see \cite{C98,DKM21,DTY23,KM08,KMZ12,KS09,DM87,DM92,DM85,MP82,MP83,MS4,MS2,MS5,Z15}
and a survey paper \cite{KM09}.  For known results of  Hamiltonicity  of  Cayley graphs, see    \cite{CG96,D83,
	GWM14,GM07,GKMM12,KM09,MR234, DM83,WM18,WM15}. Additionally, the existence of Hamilton cycle of semisymmetric graphs of order $2pq$ and $2p^3$ were recently confirmed in \cite{Chen26,DY26}, where $p$ and $q$ are primes.

 A connected simple graph is said to be of rank $r$ if  its automorphism group is a transitive group of   rank $r$ on vertices.
     Every   graph of rank 2  is a complete graph,   which contains a Hamilton cycle, and  so graphs of rank 3 are the smallest nontrivial case.
  One of our long-term goals is to show that every  such graph   contains a Hamilton cycle, except for the Petersen graph.
  It is straight to see  the the existence of   Hamilton cycles when
 the group is imprimitive, see \cite[Lemma 3.1]{DYao26}. So we only need to deal with primitive groups of rank 3.
  In this paper,  we shall deal with primitive groups of rank 3, whose socle is a classical group.
   The other remaining cases will be done in our other papers.
The following  theorem is the main result of  this paper.

\begin{theorem} \label{the:main} Let $\G$ be a graph of rank 3, whose automorphism group   is a  primitive group of rank 3, whose socle is  a classical group.
 Then  $\G $ contains a Hamilton cycle, with the exception of    the Petersen graph.
\end{theorem}

Remak that the existence of Hamilton cycles  has been confirmed   for  several graphs of rank 3, for instance,
  graphs of order $pq$ and rank 3 where $p$ and $q$ are primes, see
  \cite{DKM21}.

This paper is organized as follows: after this introductory section, some preliminaries     will be given in Section 2; the strategy of the proof of Theorem~\ref{the:main}  will be given in Section 3;  the proof can be concluded to orthogonal groups and then
the four actions    of rank 3 of orthogonal groups
 will be discussed  in Sections 4-7, separately.

\section{Preliminaries}
\label{sec:pre}
Throughout this paper, graphs are finite, simple, connected and undirected,
and groups are finite.
Furthermore, a {\em multigraph} is a generalization
of a graph in which multiple edges and loops are allowed.
Given a graph $X$, we denote the vertex set by $V(X)$ and the edge set by $E(X)$, respectively.
For two adjacent vertices $u,v \in V(X)$, we write
$u \sim v$ and denote the corresponding edge by $\{u,v\}$.

 A permutation group $G$ on a set $V$ induces
the action of $G$ on $V\times V$, whose
orbits are called {\em orbitals} of $G$, while  $O_0=\{(x,x) \colon x\in V\}$ is said to be {\it trivial}.
Every orbital $O_i$ corresponds to a {\it paired} orbital $O_i^*=\{ (x,y): (y,x)\in O_i\}$, and
$O_i$ is said to be {\it self-paired} if $O_i=O_i^*$.

For an arbitrary union $\cal{O}$ of orbitals (having empty intersection
with $O_0$), the {\em generalized orbital (di)graph} $X(V,\cal{O})$
of the action of $G$ on $V$
with respect to $\cal{O}$ is a simple (di)graph with the vertex set $V$ and
edge set $\cal{O}$.
It can be viewed as an undirected graph if and only if  $\cal{O}$ coincides with its symmetric closure,
that is, $\cal{O}$ has the property that $(x,y)\in \cal{O}$ implies $(y,x)\in \cal{O}$.
Further,  $X(V,\cal{O})$ is
said to be {\em basic orbital graph } if ${\cal O} =O_i\cup O_i^*$ for some $O_i$.

 Let $X$ be a graph with a transitive group $G$ of $\Aut(X)$. By $G_v$ and $[G:H]$ we denote the stabilizer of $v\in V(X)$ and the set of cosets of $G$ relative to $H$, respectively, where $H=G_v$.
Then there exists a 1-1 correspondence between orbitals and suborbits of $G$, that is orbits of $H$ on $V(X)$. A suborbit corresponding to a {\em self-paired}
orbital is said to be self-paired. Therefore, the (generalized) orbital (di)graph $X(V, \cal{O})$ is denoted by  $X(G,H,D)$, provided $\cal{O}$ corresponds to a union   of suborbits  (some double cosets of $H$).
 \vskip 1.5mm
Finally, we give some  known results on graphs and groups.

\begin{proposition}
	\label{pro:1}
	{\rm \cite[Theorem~1]{BJ78}}
	Every $2$-connected regular graph of order $n$ and valency at least $n/3$ contains a Hamilton cycle (every vertex-transitive graph is 2-connected, see \cite[Theorem~3]{M70}).

	\end{proposition}

\begin{proposition}
    \label{pro:2}
    {\rm \cite[Corollary~3]{A82}}
    If $X$ is a connected Cayley graph of an abelian group of order at least 3, then every edge of $X$ lies in a Hamiltonian cycle.
\end{proposition}

\begin{proposition} \label{pro:3}  {\rm \cite{KL82}} Let $\soc(F/Z(M))=M/Z(M)$ be a classical group,  acting as a primitive
rank $3$ permutation group of degree $n$, on the set $X$ of cosets of a subgroup $K/Z(M)$ of $F/Z(M)$. Then
one of the following occurs up to conjugacy under $\Aut (M/Z(M))$:
\begin{enumerate}
\item[\rm(i)] $M/Z(M)=\PSL(n,q)\le F/Z(M)\le \Aut(M/Z(M))$ and $X$ is the set of lines, $n\ge 4;$
\item[\rm(ii)]  $M$ is one of the groups
$\Sp(2m-2,q)$, $\SU(m,q)$,  $\O^{\pm}(2m,q)$,  $\O(2m-1,q)$,
for $m\ge 3$ and q a prime power,  $M\trianglelefteq F$ with $F/Z(M)\le\Aut(M/Z(M));$
\begin{enumerate}
\item[\rm(ii.1)]  $X$ is an $M$-orbit of singular (or  isotropic) points;
\item[\rm(ii.2)] $X$ is an $M$-orbit of maximal totally singular (or  isotropic) subspaces and
$M=\Sp(4, q), $ $\SU(4, q), \SU(5, q), \O^-(6, q),$ $\O^+( 8, q)$ or $\O^+(10, q).$
\item[\rm(ii.3)] $X$ is any $M$-orbit of nonsingular points and $M=\SU(m,2),$ $\O^{\pm}(2m,2)$
$\O^{\pm}(2m,3)$ or $\O(2m-1, 3)$;
\item[\rm(ii.4)]  $X$  is either orbit of nonsingular hyperplanes of $M=\O(2m-1,4)$ or $\O(2m-1,8)$
(where $G=\O(2m-1, 8)\cdot \ZZ_3$ in the latter case).
\end{enumerate}
\item[\rm(iii)] groups of small order

$M/Z(M)=\PSL(2,4)$, $|X|=10$; $M/Z(M)=\PSL(2,9)\cong A_6,$ $|X|=15$;

$M/Z(M)=\PSL(4,2)\cong A_8,$ $|X|=28$; $F/Z(M)=\PGammaL(2,8)$ $|X|=36$;

$M/Z(M)=\PSL(3,4)$, $M/Z(M)\cap (K/Z(M))=A_6$;

$M/Z(M)=\PSL(4,3)$, $M/Z(M)\cap (K/Z(M))=\PSp(4,3)$;

$M=\Sp(6,2), K=G_2(2);$  $M=\SU(3,3),$ $K\cap  M=\PSL(3,2);$

$M =\SU(3,5), K\cap M = 3\cdot A_7$; $M =\SU(4,3)$, $K\cap M=4\cdot \PSL(3,4)$;

$M=\SU(6,2), K\cap M=3\cdot \PSU(4,3)\cdot 2;$ and $M=\O(7,3), K=G_2(3)$.
\end{enumerate}
\end{proposition}

\section{Proof of   Theorem~\ref{the:main}}
\demo To prove Theorem~\ref{the:main}, let $\G$ be a  graph of rank 3, whose automorphism group   $G=\Aut(\G)$ is  a primitive  of rank 3, such that
  $\soc(G)=M/Z(M)$ is  a  classical  group.   Then $M$ is one of groups listed in Proposition~\ref{pro:3}.

  All groups  in  (3)   with small order  may be managed by Magma.

   Suppose that  $G$ is a  linear group,    symplectic group or unitary   group. Then  by \cite[Theorem 1.1]{DYao26}, either $\G$ has a Hamilton cycle, with the exception of the Petersen graph.

   All    orthogonal  groups $M$ in Proposition~\ref{pro:3} will be respectively  discussed  the following four subsections, according to its actions on  singular points, nonsingular points, maximal totally
 singular subspaces and nonsingular hyperplanes.  Remind that  when $q$ is odd,  isotropic (reps. nonisotropic) points are  singular (nonsingular).
   Thus Theorem~\ref{the:main} is proved.
\qed

\vskip 3mm

 Some basic definitions and  notations  are given in Subsections 3.1 and 3.2.
\subsection{Orthogonal spaces}
Let $V=V(n,q)$ be a  liner space of dimension $n$ over $\FF_q$.  The  quadratic form $Q$ on $V$ is a map $V\to \FF_q$ such that

(i) $Q(\ld X)=\ld ^2Q(X)$, for all $\ld \in \FF_q$ and $X\in V$; and

(ii) $f(X, Y)=Q(X+Y)-Q(X)-Q(Y)$ for all $X, Y\in V$ is a symmetric bilinear form on $V$.

Let $V$ be an orthogonal space  with  the  quadratic form $Q$,  associate with  the bilinear form $f(X,Y)$.
Let \begin{equation}\label{rad}
 \rad(f)=\{ v\in V \di f(v,w)=0,\forall w\in V \},\,  \rad(Q)=\{v\in \rad(f)\di Q(v)=0\}.\end{equation}
Then $\rad(Q)\le \rad(f)\le V$. The quadratic form $Q$ is called  {\it nondegenerate}  if $\rad(f)=0$; and  {\it nonsingular} if $\rad(Q)=0$.
Note that when $\dim(V)$ is even, $\rad(Q)=0$ if and only if $\rad(f)=0$, see \cite[Lemma 1.5.38]{BHD67}.

There is only one essential type for orthogonal space of odd dimension,  such that
 \begin{equation}\label{odd}\begin{array}{lll} &Q(X)=a^2+X_1X_2^t\quad {\rm and}\, \, &f(X,Y)=X_1Y_2^t+X_2Y_1^t, q\, \, {\rm even};\\
 &Q(X)=\frac 12 a^2+X_1X_2^t\quad {\rm and}\, \, &f(X,Y)=ab+X_1Y_2^t+X_2Y_1^t, q\, \, {\rm odd.}\end{array}
\end{equation}
for any $X=(a,X_1,X_2)$ and $Y=(b,Y_1,Y_2)$, where $X_i, Y_i\in V(m,q)$, and $^t$ is the transpose operation.

 There are two essential types  for orthogonal space of even dimension, that is plus type with
  \begin{equation}\label{plus}
  Q(X)=X_1X_2^t\,\, {\rm and}\, \, f(X,Y)=X_1Y_2^t+X_2Y_1^t.\end{equation}
 where  $X=(X_1,X_2)$ and $Y=(Y_1,Y_2)$; and minus type with
 \begin{equation}\label{minus}
\begin{array}{lll} &Q(X)=X_1X_2^t+a_1^2+a_1a_2+\mu a_2^2,\, &f(X,Y)=X_1Y_2^t+X_2Y_1^t+a_1b_2+a_2b_1,\,\, q\,  {\rm even};\\
 &Q(X)=X_1X_2^t+\frac 12(a_1^2-\th a_2^2),\, &f(X,Y)=X_1Y_2^t+X_2Y_1^t+a_1b_1-\th a_2b_2,\,\, q\,  {\rm odd}.\end{array} \end{equation}
  for a given $\mu $ such that  $x^2+x+\mu $ in irreducible,   where  $X=(X_1,X_2,a_1,a_2)$ and $Y=(Y_1,Y_2,b_1,b_2)$.

Suppose that $q$ is odd. Then $Q(X)=\frac 12f(X,X)$ so that $Q(X)=0$ if and only if $f(X,X)=0$.
In this case, isotropic points just coincide with singular points.

\subsection{Key notations}
Throughout the paper,  let $\FF_q^*=\lg \th  \rg$ and  $\FF_{q^n}^*=\lg \ld \rg$, where $n\ge 2$,  $\ld ^n=\sum_{i=0}^{n-1}a_i\ld ^i$. Then $\overline {\ld }(x)=\ld x$,  $\forall x\in \FF_{q^n}$  gives a linear transformation of $\FF_{q^n}$ over $\FF_q$.
Moreover, the matrix $A$ of $\overline{\ld }$  under the base: $1, \ld, \cdots, \ld^{n-1}$ is
\begin{equation}\label{A1} A={\left(\begin{array}{cccccc}0 &1 &0 &\cdots &0 &0\\ 0 &0 &1 &\cdots &0 & 0 \\  &&&\cdots&&\\0 &0 &0  &\cdots &0&1 \\   a_0 &a_1 &a_2 &\cdots
&a_{n-2} &a_{n-1} \end{array}\right)}\end{equation}
and then
\begin{equation}\label{A2} A^{-1}={\left(\begin{array}{ccccc} -a_1a_0^{-1} &-a_2a_0^{-1} &\cdots &-a_{n-1}a_0^{-1} & a_0^{-1}\\ 1 &0 &\cdots &0 &0\\ 0 &1 &\cdots &0 &0\\ &&\cdots &&\\0 &0 &\cdots &1 &0\end{array}\right)}.\end{equation}

 By $\Del_s$ and $\Del_l$ we denote the shorter and longer one, respectively, where $|\Del_s|$  may be equal to  $|\Del_l|.$ Then  $|\Del_l|\ge \frac{|V|}3.$
  By Proposition~\ref{pro:1}  the corresponding basic orbital graph $\G(V, \Del_l)$ has a Hamilton cycle.  So we one need to prove the theorem for  $\G(V, \Del_s)$. Later we always denote $V$ by $\O$. So
 set
\begin{equation} \label{orbital} \G=\G(\O, \Del_s) .  \end{equation}

 For any H-path leading $\a$ to $\b$  of a certain graph,  we can simply denote that by $\a \to \b$. From now on,   use the abbreviation   H-cycle (resp. H-path) for  Hamilton cycle (resp. Hamilton path).
Finally  by $A^\t$  we denote   the inverse transpose of $A$.

 \section{Actions on  singular points}
By Proposition~\ref{pro:3}, we consider    the action of $M$ on   singular points by dealing with  three  kinds of   orthogonal space $V(n,q)$   separately, refer Eq(\ref{odd})-(\ref{minus}).

\subsection{$V$ is of odd dimension}

\begin{lem}\label{os} Let $V=V(n,q)$ be the orthogonal space, where $n=2m+1\ge 5$.
Let $M=\O(n, q)$ acting on the  set $\O$ of  singular points. Then  $\G(\O, \Del_s)$  has a H-cycle.
\end{lem}
\demo  Since $\P\O(5,q)\cong \PSp(4,q)$, which was discussed in  \cite[Lemma 5.1.(2)]{DYao26}, we assume   $m\ge 3$ and deal with two cases when $q$ is odd and even, separately.
\vskip 2mm
{\bf (1)} $q$ is odd. \vskip 2mm

By Eq(\ref{odd}), set the  bilinear form $f$ as
  $$f(X,Y)=ab+X_1Y_2^t+X_2Y_1^t,$$
where  $X=(a,X_1,X_2),  Y=(b,Y_1,Y_2), \, X_i,\, Y_i \in V(m,q).$
Let $\O$ be the set of singular points so that   $|\O|=\frac{q^{n-1}-1}{q-1}$. Take  $\a \in \O$.
 Then $\D_s=\{\b \in \O \di f(\a,\b)=0\}$ where $|\D_s|=\frac{q(q^{n-3}-1)}{q-1}$.
Let $\G=\G (\O, \D_s)$, see \ref{orbital}. Then $\a, \b\in E(\G)$ if and only if $f(\a, \b)=0$.

  Let $A$ and  $A^{-1}$ be the matrices   given in Eq(\ref{A1}, \ref{A2}) by replacing $n$ by $m$. Set
  $\rho=\diag((1), A, A^\t)$, a block diagonal matrix.  Then $\rho \in \SO(n,q)$ and $\rho(X)=(a,X_1A,X_2A^\t)$, for any $X=(a, X_1, X_2)$. Set

\begin{equation} \label{odd2}
\begin{array}{llll}
 &r&=\frac{q^m-1}{q-1},\quad 0_m=(0, \cdots , 0 ), &X_0=(0,1,0,\cdots,0)\in V(m,q), \\
&Y&=(x_1, x_3, \cdots , x_m)\in V(m-1,q),  &X_{e, Y}=(x_1, e, x_3, \cdots x_m)\in V(m, q).
\end{array}
\end{equation}
Set
$$\begin{array}{llll}
B_1&=&\{\b_{1,i}\di 0\le i \le r-1\}, \quad  &\b_{1,i}=\lg(0,X_0A^i,0_m)\rg ;\\
B_2&=&\{\b_{2,i}\di 0\le i \le r-1\}, \quad  &\b_{2,i}=\lg(0,0_m,X_0(A^i)^\t)\rg ;\\
C_{Y}&=&\{\d_{Y, k}\di 0\le k \le \frac {q^m-1}{2}-1\}, &\d_{Y, k}=\lg(0,X_0A^k,X_{0,Y}(A^k)^\t)\rg , \\
&&& Y=\lg \th^2\rg Y_1, Y_1\in V(m-1,q)\setminus \{0_{m-1}\};\\
D_{Y}&=&\{\g_{Y, k}\di 0\le k \le q^m-2\}, &\g_{Y, k}=\lg(1,X_0A^k,X_{-\frac{1}{2}, Y}{(A^k)^\t})\rg ,\\
&&& Y\in V(m-1,q).
\end{array} $$
Check that $\O$ is a disjoint union of the above  $\lg\rho\rg$-orbits, that is
$$\O=B_1\cup B_2\cup (\cup _{Y}C_Y)\cup (\cup _{Y}D_Y),$$
where  the respective vectors $Y$ for  $C_Y$ and $D_Y$ are given as above.
Moreover, $\G(B_1)$, $\G(B_2)$, $\G(\{\d_{ Y, k}\di Y\})$ and $\G(\{\g_{Y, k}\di Y\})$ are complete graphs,  for any given $k$.
\vskip 2mm
{\it Step i:}  Set
$$Y_1=(1, 0,\cdots,0),\,\, Y_2=(0,1,\cdots,0), \, \,  Y_3=(1,1,0,\cdots,0)\in V(m-1, q).$$  Since $X_0A^{-1}X_{0, Y_2}^t+X_0AX_{0, Y_1}^t=0$, we get $\d_{Y_1, 0}\thicksim \d_{Y_2, 1}=\d_{Y_2, 0}^{\orho}$.
In the complete graph $\G(\{\d_{Y, 0}\di Y\})$, pick a H-path
$$U_0:\quad \d_{Y_2, 0}, \cdots , \d_{Y_3, 0}, \d_{Y_1, 0},$$
and then get a H-path $U$ of  $\G(\cup_Y C_{Y})$ (see the left part of Figure 1.):
$$U:\quad \d_{Y_1, 0},U_0^{\orho},U_0^{\orho^2},\cdots,U_0^{\orho^{\frac {q-1}{2}r-1}},\d_{Y_2, 0}, \cdots, \d_{Y_3, 0}.$$

{\it Step ii:} Set
$$Y_4=(-\frac{1}{2},0,0,\cdots,0)\, \, {\rm  and}\, \,  Y_5=(0,-\frac{1}{2},0,\cdots,0)\in V(m-1,q).$$
Then $\g_{Y_5, 0}\thicksim\g_{Y_4, 1}=\g_{Y_4, 0}^\rho$. In the complete graph $\G(\{\g_{Y, 0}\di Y\})$,
pick a H-path
 $$W_0:\quad \g_{Y_4, 0}, \g_{Y_3, 0},  \cdots ,    \g_{Y_5, 0},$$
  and then get a H-cycle $W$ of $\G(\cup_{Y}D_{Y})$  (see the right part of Figure 1.):
$$W:\quad \g_{Y_3, 0},  \cdots ,    \g_{Y_5, 0}, W_0^{\rho},W_0^{\rho^2},\cdots,W_0^{\rho^{(q-1)r-1}},\g_{Y_4, 0}.$$

{\it Step iii:} Take a H-path $\b_{2,1}\to\b_{2, r-1}$ of $\G(B_2)$ and a H-path $\b_{1,1}\to\b_{1,0}$ of  $\G(B_1)$.
Finally,  we get a H-cycle $C$ for $\G$:
$$C:\quad U, \b_{2,1}\to\b_{2,r-1}, W, \b_{1,1}\to\b_{1,0},\d_{Y_1, 0}.$$

\begin{figure}[htbp]
 \vspace{-2mm} \centering
  \begin{tikzpicture}
  [
    scale=0.75, transform shape,
    p/.style={draw, circle, fill=white, inner sep=0pt, minimum size=2.6mm},
    d/.style={fill=white, inner sep=2pt},
    box/.style={draw, rounded corners=6pt},
    lab/.style={anchor=north, yshift=3pt, font=\footnotesize},
    tlab/.style={anchor=south, font=\small},
    midarrow/.style={
      postaction=decorate,
      decoration={markings, mark=at position 0.5 with {\arrow{angle 45}}}
    },
    remidarrow/.style={
      postaction=decorate,
      decoration={markings, mark=at position 0.5 with {\arrowreversed{angle 45}}}
    }
  ]
  \node[p] (a10) at (0,0) {};
  \node[d] (d10) at (2,0) {$\cdots$};
  \node[p] (a20) at (4,0) {};
  \node[p] (a30) at (6,0) {};
  \node[p] (c1) at (8,0) {};
  \node[d] (d20) at (8.5,0) {$\cdots$};
  \node[p] (c2) at (9,0) {};
  \node[p] (b10) at (11,0) {};
  \node[p] (b20) at (13,0) {};
  \node[d] (d30) at (15,0) {$\cdots$};
  \node[p] (b30) at (17,0) {};
  \draw[midarrow] (a10.east) -- (d10.west);
  \draw[midarrow] (d10.east) -- (a20.west);
  \draw[midarrow] (b20.east) -- (d30.west);
  \draw[midarrow] (d30.east) -- (b30.west);
  \node[p] (a11) at (0,-1) {};
  \node[d] (d11) at (2,-1) {$\cdots$};
  \node[p] (a21) at (4,-1) {};
  \node[p] (a31) at (6,-1) {};
  \node[p] (d1) at (8,-1.2) {};
  \node[d] (d21) at (8.5,-1.2) {$\cdots$};
  \node[p] (d2) at (9,-1.2) {};
  \node[p] (b11) at (11,-1) {};
  \node[p] (b21) at (13,-1) {};
  \node[d] (d31) at (15,-1) {$\cdots$};
  \node[p] (b31) at (17,-1) {};
  \draw[midarrow] (a11.east) -- (d11.west);
  \draw[midarrow] (d11.east) -- (a21.west);
  \draw[midarrow] (a21.east) -- (a31.west);
  \draw[midarrow] (b11.east) -- (b21.west);
  \draw[midarrow] (b21.east) -- (d31.west);
  \draw[midarrow] (d31.east) -- (b31.west);
  % ...
  \node[d] at (0,-2) {$\vdots$};
  \node[d] at (2,-2) {$\vdots$};
  \node[d] at (4,-2) {$\vdots$};
  \node[d] at (6,-2) {$\vdots$};

  \node[d] at (11,-2) {$\vdots$};
  \node[d] at (13,-2) {$\vdots$};
  \node[d] at (15,-2) {$\vdots$};
  \node[d] at (17,-2) {$\vdots$};
  \node[p] (a12) at (0,-3) {};
  \node[d] (d12) at (2,-3) {$\cdots$};
  \node[p] (a22) at (4,-3) {};
  \node[p] (a32) at (6,-3) {};

  \node[p] (b12) at (11,-3) {};
  \node[p] (b22) at (13,-3) {};
  \node[d] (d32) at (15,-3) {$\cdots$};
  \node[p] (b32) at (17,-3) {};
  \draw[midarrow] (a12.east) -- (d12.west);
  \draw[midarrow] (d12.east) -- (a22.west);
  \draw[midarrow] (a22.east) -- (a32.west);
  \draw[midarrow] (b12.east) -- (b22.west);
  \draw[midarrow] (b22.east) -- (d32.west);
  \draw[midarrow] (d32.east) -- (b32.west);

  \draw[midarrow] (a30.west) -- (a11.east);
  \draw[midarrow] (a31.west) -- ($(a31)!0.36!(a12)$);
  \draw[midarrow] ($(a31)!0.64!(a12)$) -- (a12.east);
  \draw[midarrow, rounded corners=6pt]
  (a32.south) -- (6,-3.9) -- (-0.9,-3.9) -- (-0.9,0) -- (a10.west);

  \draw[midarrow] (b30.west) -- (b11.east);
  \draw[midarrow] (b31.west) -- ($(b31)!0.36!(b12)$);
  \draw[midarrow] ($(b31)!0.64!(b12)$) -- (b12.east);
  \draw[midarrow, rounded corners=6pt]
  (b32.south) -- (17,-3.9) -- (10.1,-3.9) -- (10.1,0) -- (b10.west);

  \draw[midarrow, rounded corners=6pt]
  (a20.north) -- (4,0.7) -- (8,0.7) -- (c1.north);
  \draw[midarrow, rounded corners=6pt]
  (c2.north) -- (9,0.7) -- (13,0.7) -- (b20.north);

  \draw[midarrow] (b10.west) -- (d2.east);
  \draw[midarrow] (d1.west) -- (a30.east);

  \draw[box] (-0.7,0.5) rectangle (0.7,-3.7);
  \draw[box] (1.3,0.5) rectangle (2.7,-3.7);
  \draw[box] (3.3,0.5) rectangle (4.7,-3.7);
  \draw[box] (5.3,0.5) rectangle (6.7,-3.7);

  \draw[box] (10.3,0.5) rectangle (11.7,-3.7);
  \draw[box] (12.3,0.5) rectangle (13.7,-3.7);
  \draw[box] (14.3,0.5) rectangle (15.7,-3.7);
  \draw[box] (16.3,0.5) rectangle (17.7,-3.7);

  \draw[box] (7.3,0.5) rectangle (9.7,-0.5);
  \draw[box] (7.3,-0.7) rectangle (9.7,-1.7);
  \node[lab] at (a10.south) {$\d_{Y_2,0}$};
  \node[lab] at (a20.south) {$\d_{Y_3,0}$};
  \node[lab] at (a30.south) {$\d_{Y_1,0}$};
  \node[lab] at (a11.south) {$\d_{Y_2,1}$};
  \node[lab] at (a21.south) {$\d_{Y_3,1}$};
  \node[lab] at (a31.south) {$\d_{Y_1,1}$};
  \node[lab] at (a12.south) {$\d_{Y_2,\frac{q^m-1}{2}-1}$};
  \node[lab] at (a22.south) {$\d_{Y_3,\frac{q^m-1}{2}-1}$};
  \node[lab] at (a32.south) {$\d_{Y_1,\frac{q^m-1}{2}-1}$};

  \node[lab] at (b10.south) {$\d_{Y_4,0}$};
  \node[lab] at (b20.south) {$\d_{Y_3,0}$};
  \node[lab] at (b30.south) {$\d_{Y_5,0}$};
  \node[lab] at (b11.south) {$\d_{Y_4,1}$};
  \node[lab] at (b21.south) {$\d_{Y_3,1}$};
  \node[lab] at (b31.south) {$\d_{Y_5,1}$};
  \node[lab] at (b12.south) {$\d_{Y_4,q^m-2}$};
  \node[lab] at (b22.south) {$\d_{Y_3,q^m-2}$};
  \node[lab] at (b32.south) {$\d_{Y_5,q^m-2}$};

  \node[lab] at (c1.south) {$\b_{2,1}$};
  \node[lab] at (c2.south) {$\b_{r-1}$};
  \node[lab] at (d1.south) {$\b_{1,0}$};
  \node[lab] at (d2.south) {$\b_{1,1}$};

  \node[tlab] at (0,0.6) {$C_{Y_2}$};
  \node[tlab] at (2,0.6) {$\cdots$};
  \node[tlab] at (4,0.6) {$C_{Y_3}$};
  \node[tlab] at (6,0.6) {$C_{Y_1}$};

  \node[tlab] at (11,0.6) {$D_{Y_4}$};
  \node[tlab] at (13,0.6) {$D_{Y_3}$};
  \node[tlab] at (15,0.6) {$\cdots$};
  \node[tlab] at (17,0.6) {$C_{Y_5}$};

  \node[tlab] at (8,0.6) {$B_2$};
  \node[tlab] at (8,-2.3) {$B_2$};
  \end{tikzpicture}
  \caption{{\footnotesize H-cycle C of the graph $\G$ in Lemma 4.1.{\bf (1)}}}
\end{figure}

{\bf (2)} $q$ is even.
\vskip 2mm
By Eq(\ref{odd}), set the  quadratic form $Q$ and bilinear form $f$ as
  $$Q(X)=a^2+X_1X_2^t \, \,  \, {\rm and}\,\, f(X,Y)=X_1Y_2^t+X_2Y_1^t,$$
  where  $ X=(a,X_1,X_2), \, Y=(b,Y_1,Y_2), \, X_i, Y_i\in V(m,q).$
  Let $\O$ be the set of  singular points $\lg v\rg $ (that is $Q(v)=0$),  where $|\O|=\frac{q^{n-1}-1}{q-1}$.
Take   $\a \in \O$. Then $\D_s=\{\b \in \O \di f(\a,\b)=0\}$, where $|\D_s|=\frac{q(q^{n-3}-1)}{q-1}$.

 Let $A, \rho,  r, 0_m, X_0$ and $X_{e,Y}$  have the same meanings in Eq(\ref{odd2}).
Set
$$\begin{array}{llll}
B_1&=&\{\b_{1,i}\di 0\le i \le r-1\}, \quad  &\b_{1,i}=\lg(0,X_0A^i,0_m)\rg ;\\
B_2&=&\{\b_{2,i}\di 0\le i \le r-1\}, \quad  &\b_{2,i}=\lg(0,0_m,X_0(A^i)^\t)\rg ;\\
C_{Y}&=&\{\d_{Y, k}\di 0\le k \le q^m-2\}, &\d_{Y, k}=\lg(0,X_0A^k,X_{0,Y}(A^k)^\t)\rg , \\
&&& Y=\lg \th\rg Y_1, Y_1\in V(m-1,q)\setminus \{0_{m-1}\};\\
D_{Y}&=&\{\g_{Y, k}\di 0\le k \le q^m-2\}, &\g_{Y, k}=\lg(1,X_0A^k,X_{1, Y}{(A^k)^\t})\rg ,\\
&&& Y\in V(m-1,q).
\end{array} $$
Now $\O$ is a disjoint union of the above  $\lg\rho\rg$-orbits:
$$\O=B_1\cup B_2\cup (\cup _{Y}C_Y)\cup (\cup _{Y}D_Y),$$
where  the respective vectors $Y$ for  $C_Y$ and $D_Y$ are given as above.
 Set $$Y_4=Y_1=(1, 0, \cdots,  0),\, \, Y_5=Y_2=(0, 1, \cdots,  0)\, \, {\rm and}\, \,   Y_3=(1, 1, \cdots,  0).$$
Then the  remaining arguments are completely same as that in   {\bf (1)}. \qed

\subsection{$V$ is of plus type}

\begin{lem}\label{esp}
Let $V=(2m,q)$ be an orthogonal space of   plus type, where $m\ge 3$.  Let $M=\O(2m, q)$ acting on the set  $\O$ of  singular  points.
Then  $\G(\O, \Del_s)$  has a H-cycle.
\end{lem}
\demo
 {\bf (1)}  $q$ is odd. \vskip 2mm

By Eq(\ref{plus}), set
  $$f(X,Y)=X_1Y_2^t+X_2Y_1^t,$$
  where $X=(X_1,X_2),Y=(Y_1,Y_2)$,  $X_i,Y_i \in V(m,q)$.  Let $\O$ be the set of   singular  points so that   $|\O|=\frac{(q^m-1)(q^{m-1}+1)}{q-1}$.
 Take $\a \in \O$. Then $\D_s=\{\b \in \O \di f(\a,\b)=0\}$, where $|\D_s|=\frac {q(q^{m-1}-1)(q^{m-2}+1)}{q-1}$. Let $\G=\G(\O, \D_s)$.

Let $A$ be the matrix given in Eq(\ref{A1}), replacing $n$ by $m$. Let $\rho=\diag (A, A^\t)$. Then $\rho \in \SO(n,q)$ and $\rho(X)=(X_1A,X_2A^\t)$.
Adapt the notations in Eq(\ref{odd2}).  Set
$$\begin{array}{llll}
B_1&=&\{\b_{1,i}\di 0\le i \le r-1\}, \quad  &\b_{1,i}=\lg(X_0A^i,0_m)\rg ;\\
B_2&=&\{\b_{2,i}\di 0\le i \le r-1\}, \quad  &\b_{2,i}=\lg(0_m,X_0(A^i)^\t)\rg ;\\
C_{Y}&=&\{\d_{Y, k}\di 0\le k \le \frac {q^m-1}{2}-1\}, &\d_{Y, k}=\lg(X_0A^k,X_{0,Y}(A^k)^\t)\rg , \\
&&& Y=\lg \th^2\rg Y_1, Y_1\in V(m-1,q)\setminus \{0_{m-1}\}.\\
\end{array} $$
Check that $\O$ is a disjoint union of the above $\lg\rho\rg$-orbits, that is
$$\O=B_1\cup B_2\cup (\cup _{Y}C_Y).$$
Moreover,  $\G(B_1)$, $\G(B_2)$ and $\G(\{\d_{Y,k} \di Y\})$ are complete graphs, for any given $k$.
Set
 $$Y_1=(1, 0,\cdots,0), \, \, Y_2=(0,1,\cdots,0)\, \, {\rm  and}\, \,   Y_3=(1,1,0,\cdots,0)\in V(m-1, q).$$
   Check that  $\d_{Y_1, 0}\thicksim \d_{Y_2, 1}=\d_{Y_2, 0}^{\orho}$.
In the complete graph $\G(\{\d_{Y, 0}\di Y\})$, pick a H-path:
$$U_0:\quad \d_{Y_2, 0}, \cdots ,\d_{Y_3,0},  \d_{Y_1, 0},$$ and then get a H-cycle in $\G(\cup_Y C_{Y})$:
$$U:\quad \d_{Y_1, 0}, U_0^{\orho},U_0^{\orho^2},\cdots,U_0^{\orho^{\frac {q-1}{2}r-1}},\d_{Y_2, 0}, \cdots ,  \d_{Y_3, 0}.$$
Take  a H-path $\b_{2,1}\to\b_{2,0}$ in $\G(B_2)$ and a H-path $\b_{1,1}\to\b_{1,0}$ in $\G(B_1)$.
Then  we get a H-cycle $C$ for $\G$:
$$C:\quad U, \b_{2,1}\to\b_{2,0},\b_{1,1}\to\b_{1,0},\d_{Y_1, 0}.$$

{\bf (2)} $q$ is even. \vskip 2mm

By Eq(\ref{plus}),  set
$$Q(X)=X_1X_2^t\, \, \, {\rm and}\, \,\,  f(X,Y)=X_1Y_2^t+X_2Y_1^t, $$
 where  $X=(X_1,X_2)$ and $Y=(Y_1,Y_2)$.
  Let $\O$ be the set of   singular  points  so that    $|\O|=\frac{(q^m-1)(q^{m-1}+1)}{q-1}$.
Take  $\a\in \O$. If $q=2$, then $\D_s=\{\b \in \O \di f(\a,\b)\ne0\}$, where $|\D_s|\ge |\O|/3$.
 So assume  $q\ge3$ so that   $\D_s=\{\b \in \O \di f(\a,\b)=0\}$, where $|\D_s|=\frac {q(q^{m-1}-1)(q^{m-2}+1)}{q-1}$.

Let $A,  r, 0_m, X_0$ and $X_{e,Y}$  have the same meanings as in (1) (see Eq{(\ref{odd2})}).  Then  $\rho \in \SO(n,q)$ and $\rho(X)=(X_1A,X_2A^\t)$.
Set
$$\begin{array}{llll}
B_1&=&\{\b_{1,i}\di 0\le i \le r-1\}, \quad  &\b_{1,i}=\lg(X_0A^i,0_m)\rg ;\\
B_2&=&\{\b_{2,i}\di 0\le i \le r-1\}, \quad  &\b_{2,i}=\lg(0_m,X_0(A^i)^\t)\rg ;\\
C_{Y}&=&\{\d_{Y, k}\di 0\le k \le q^m-2\}, &\d_{Y, k}=\lg(X_0A^k,X_{0,Y}(A^k)^\t)\rg , \\
&&& Y=\lg \th\rg Y_1, Y_1\in V(m-1,q)\setminus \{0_{m-1}\}.\\
\end{array} $$
Check that $\O$ is a disjoint union of the above  $\lg\rho\rg$-orbits, that is $\O=B_1\cup B_2\cup (\cup _{Y}C_Y).$
The remaining  argument just follows {\bf (1)}.
\qed
\subsection{$V$ is of minus type}
\begin{lem}\label{esm}
Let $V(2m,q)$ be an orthogonal space of  minus type, where $m\ge 3$.  Let $M=\O(2m, q)$ acting on the   set of $\O$ of  singular points.
Then  $\G(\O, \Del_s)$  has a H-cycle.
\end{lem}
\demo Since  $\P\O^{-}(6,q)\cong \PSU(4,q)$, see \cite[Lemma 6.1]{DYao26}, we  assume $m\ge 4$.
\vskip 1.5mm
{\bf (1)} $q$ is odd. \vskip 1.5mm

By Eq(\ref{minus}), set
$$f(X,Y)=X_1Y_2^t+X_2Y_1^t+a_1b_1-\th a_2b_2,$$
where  $X=(X_1,X_2,a_1,a_2),Y=(Y_1,Y_2,b_1,b_2)$,   $X_i,Y_i \in V(m-1,q)$, $a_i,b_i\in\FF_q$ and   $\FF_q^*=\lg \th \rg $.
Let $\O$ be the set of   singular  points so that   $|\O|=\frac{(q^m+1)(q^{m-1}-1)}{q-1}$.
 Take    $\a \in \O$. Then $\D_s=\{\b \in \O \di f(\a,\b)=0\}$, where $|\D_s|=\frac {q(q^{m-1}+1)(q^{m-2}-1)}{q-1}$.

Let $A$ be the matrix given in Eq(\ref{A1}), replacing $n$ by $m-1$.
Let $\rho=\diag (A, A^\t, 1, 1)$.  Then $\rho \in  \SO(n,q)$ and $\rho(X)=(X_1A,X_2A^\t, a_1, a_2)$.  Use the notations in Eq(\ref{odd2}), set
$$\begin{array}{llll}
B_1&=&\{\b_{1,i}\di 0\le i \le r-1\}, \quad  &\b_{1,i}=\lg(X_0A^i,0_{m-1},0,0)\rg ;\\
B_2&=&\{\b_{2,i}\di 0\le i \le r-1\}, \quad  &\b_{2,i}=\lg(0_{m-1},X_0(A^i)^\t,0,0)\rg ;\\
C_{Y}&=&\{\d_{Y, k}\di 0\le k \le \frac {q^{m-1}-1}{2}-1\}, &\d_{Y, k}=\lg(X_0A^k,X_{0,Y}(A^k)^\t,0,0)\rg , \\
&&& Y=\lg \th^2\rg Y_1, Y_1\in V(m-2,q)\setminus \{0_{m-2}\};\\
D_{Y}&=&\{\g_{Y, k}\di 0\le k \le q^{m-1}-2\}, &\g_{Y, k}=\lg(X_0A^k,X_{\frac{1}{2}\th , Y}{(A^k)^\t}, 0, 1)\rg ,\\
&&& Y\in V(m-2,q);\\
E_{b,Y}&=&\{\e_{b, Y; k}\di 0\le k \le q^{m-1}-2\}, &\e_{b, Y; k}=\lg(X_0A^k,X_{-\frac{1}{2}(1-b^2\th), Y}{(A^k)^\t}, 1, b)\rg ,\\
&&& Y\in V(m-2,q), b\in \FF_q.
\end{array} $$
Check that $\O$ is a disjoint union of the above $\lg\rho\rg$-orbits, that is
$$\O=B_1\cup B_2\cup (\cup _{Y}C_Y)\cup (\cup _{Y}D_Y) \cup (\cup _{b, Y}E_{b,Y}),\vspace{-2mm}$$
where  the respective vectors $Y$ for  $C_Y$, $D_Y$ and   $E_{b,Y}$ are given as above, and $b\in \FF_q$.
Moreover,  $\G(B_1)$, $\G(B_2)$, $\G(\{\d_{Y, k}\di Y\})$, $\G(\{\g_{Y,k}\di Y\})$ for any given $k$ and $\G(\e_{b,Y; k}\di Y)$ for any given $k, b$, are complete graphs.
Then we find H-cycle of $\G$ by the following five steps.

\vskip 2mm
{\it Step i:}  Let
$$Y_1=(1,0,0,\cdots,0),\,\,  Y_2=(0,1,0,\cdots,0), \, \, Y_3=(1,1,0,\cdots,0)\in V(m-2,q).$$
 Then  $\d_{Y_1,0}\thicksim\d_{Y_2,1}$. In the complete graph $\G(\{\d_{Y,0}\di Y\})$, pick a H-path:
 $$U_0:\quad \d_{Y_2,0}, \cdots ,  \d_{Y_3,0}, \d_{Y_1,0},$$  and then get a H-path $U$  in $\G(\cup_{Y} C_Y)$ (see the left part of Figure 3):
 $$U:\quad \d_{Y_1,0} ,U_0^{\orho},U_0^{\orho^2},\cdots,U_0^{\orho^{\frac {q-1}{2}r-1}}, \d_{Y_2,0}, \cdots ,  \d_{Y_3,0}.$$

{\it Step ii:} Let $$Y_4=(\frac{1}{2}\th,0,0,\cdots,0)\,\, {\rm  and}\,\,  Y_5=(0,\frac{1}{2}\th,0,\cdots,0)\in V(m-2,q).$$ Then  $\g_{Y_5, 0}\thicksim\g_{Y_4, 1}=\g_{Y_4, 0}^\rho$. In the complete graph $\G(\{\g_{Y, 0}\di Y\})$,
pick a H-path:
 $$W_0:\quad \g_{Y_4, 0}, \g_{Y_3, 0}, \cdots ,  \g_{Y_5, 0},$$ and then get a H-cycle  $W$ in $\G(\cup_{Y}D_{Y})$  (see the left part of Figure 2):
$$W:\quad \g_{Y_3, 0}, \cdots ,  \g_{Y_5, 0}, W_0^{\rho},W_0^{\rho^2},\cdots,W_0^{\rho^{(q-1)r-1}}, \g_{Y_4, 0}.$$

{\it Step iii:} Write $\FF_q=\{ b_i\di 0\le i\le q-1\}$. For every $b_i$, let $$Y_{6}(b_i)=(-\frac{1}{2}(1-b_i^2\th),0,0,\cdots,0),\,\, Y_7(b_i)=(0,-\frac{1}{2}(1-b_i^2\th),0,\cdots,0)\in V(m-2,q).$$ Then  we get $\e_{b_i,Y_{7}(b_i); 0}\thicksim \e_{b_i,Y_{6}(b_i);1}$.
In the complete graph $\G(\{\e_{b_i,Y;0}\di Y\})$, pick a H-path:
$$W_0(b_i):\quad  \e_{b_i,Y_{6}(b_i);0}, \e_{b_i,Y_3;0}, \cdots , \e_{b_i,Y_{7}(b_i);0},$$
 and then get a H-path in $\G(\cup_{Y}E_{b_i,Y})$ (see the right part of Figure 2):
$$W(b_i):\quad \e_{b_i,Y_3;0}, \cdots , \e_{b_i,Y_{7}(b_i);0}, W_{0}(b_i)^{\rho},W_{0}(b_i)^{\rho^2},\cdots,W_{0}(b_i)^{\rho^{(q-1)r-1}},\e_{b_i,Y_{6}(b_i); 0}.$$
We may  rechoose a H-path   $W(b_i)^{\rho^{2(i+1)}}$
in $\G(\cup_{Y}E_{b_i,Y})$ as shown in  the right part of Figure 3. Then attach $W$ and these  $W(b_i)^{\rho^{2(i+1)}}$ ($0\le i\le q-1$) by some vertices $\b_{2,1}, \b_{2,3}, \cdots ,$ $\b_{2, 2q-1}$. Thus
we get a path $R$ (see the right part of Figure 3):
$$R:\quad W,\b_{2,1}, W(b_0)^{\rho^2}, \b_{2,3}, W(b_1)^{\rho^4},
\cdots,W(b_{q-2})^{\rho^{2(q-1)}},\b_{2,2q-1}, W(b_{q-1})^{\rho^{2q}}.$$

{\it Step iv:} Let $B_2'=B_2\setminus \{\{\b_{2,1}, \b_{2,3}, \cdots ,$ $ \b_{2, 2q-1}\}\cup \{\b_{2, r-1}\}\}$. Take a H-path $\b_{2,2q+1}\to\b_{2,0}$ of the complete graph $\G(B_2')$;
and a H-path     of the complete graph $\b_{1,r-1}\to\b_{1,0}$ in $\G(B_1)$.
\vskip 2mm
{\it Step v:} Attaching the above paths from Steps i-iv,  we finally get a H-cycle $C$ for $\G$ (see Figure 3):
$$C:\quad U, \b_{2,r-1},R,\b_{2,2q+1}\to\b_{2,0},\b_{1,r-1}\to\b_{1,0},\d_{Y_1,0}.$$

\begin{figure}[htbp]
 \vspace{-2mm}  \centering
  \begin{tikzpicture}
  [
    scale=0.75, transform shape,
    p/.style={draw, circle, fill=white, inner sep=0pt, minimum size=2.6mm},
    d/.style={fill=white, inner sep=2pt},
    box/.style={draw, rounded corners=6pt},
    lab/.style={anchor=north, yshift=3pt, font=\footnotesize},
    tlab/.style={anchor=south, font=\small},
    midarrow/.style={
      postaction=decorate,
      decoration={markings, mark=at position 0.5 with {\arrow{angle 45}}}
    },
    remidarrow/.style={
      postaction=decorate,
      decoration={markings, mark=at position 0.5 with {\arrowreversed{angle 45}}}
    }
  ]
    \node[p] (a10)  at (0, 0)  {};
    \node[p] (a20)  at (2, 0)  {};
    \node[d] (d10)  at (4, 0)  {$\cdots$};
    \node[p] (a30)  at (6, 0)  {};
    \node[p] (b10)  at (9, 0)  {};
    \node[p] (b20)  at (11, 0)  {};
    \node[d] (d20)  at (13, 0)  {$\cdots$};
    \node[p] (b30)  at (15, 0)  {};
    \draw[midarrow] (a20.east) -- (d10.west);
    \draw[midarrow] (d10.east) -- (a30.west);
    \draw[midarrow] (b20.east) -- (d20.west);
    \draw[midarrow] (d20.east) -- (b30.west);
    \node[p] (a11)  at (0, -1)  {};
    \node[p] (a21)  at (2, -1)  {};
    \node[d] (d11)  at (4, -1)  {$\cdots$};
    \node[p] (a31)  at (6, -1)  {};
    \node[p] (b11)  at (9, -1)  {};
    \node[p] (b21)  at (11, -1)  {};
    \node[d] (d21)  at (13, -1)  {$\cdots$};
    \node[p] (b31)  at (15, -1)  {};
    \draw[midarrow] (a11.east) -- (a21.west);
    \draw[midarrow] (a21.east) -- (d11.west);
    \draw[midarrow] (d11.east) -- (a31.west);
    \draw[midarrow] (b11.east) -- (b21.west);
    \draw[midarrow] (b21.east) -- (d21.west);
    \draw[midarrow] (d21.east) -- (b31.west);
    \node[d] (v1)  at (0, -2)  {$\vdots$};
    \node[d] (v2)  at (2, -2)  {$\vdots$};
    \node[d] (v3)  at (4, -2)  {$\vdots$};
    \node[d] (v4)  at (6, -2)  {$\vdots$};
    \node[d] (v5)  at (9, -2) {$\vdots$};
    \node[d] (v6)  at (11, -2) {$\vdots$};
    \node[d] (v7)  at (13, -2) {$\vdots$};
    \node[d] (v8)  at (15, -2) {$\vdots$};
    \node[p] (a12)  at (0, -3)  {};
    \node[p] (a22)  at (2, -3)  {};
    \node[d] (d12)  at (4, -3)  {$\cdots$};
    \node[p] (a32)  at (6, -3)  {};
    \node[p] (b12)  at (9, -3)  {};
    \node[p] (b22)  at (11, -3)  {};
    \node[d] (d22)  at (13, -3)  {$\cdots$};
    \node[p] (b32)  at (15, -3)  {};
    \draw[midarrow] (a12.east) -- (a22.west);
    \draw[midarrow] (a22.east) -- (d12.west);
    \draw[midarrow] (d12.east) -- (a32.west);
    \draw[midarrow] (b12.east) -- (b22.west);
    \draw[midarrow] (b22.east) -- (d22.west);
    \draw[midarrow] (d22.east) -- (b32.west);

    \draw[midarrow] (a30.west) -- (a11.east);
    \draw[midarrow] (a31.west) -- ($(a31)!0.36!(a12)$);
    \draw[midarrow] ($(a31)!0.64!(a12)$) -- (a12.east);
    \draw[midarrow, rounded corners=6pt]
    (a32.south) -- (6,-3.7) -- (-1,-3.7) -- (-1,0) -- (a10.west);
    \draw[midarrow] (b30.west) -- (b11.east);
    \draw[midarrow] (b31.west) -- ($(b31)!0.36!(b12)$);
    \draw[midarrow] ($(b31)!0.64!(b12)$) -- (b12.east);
    \draw[midarrow, rounded corners=6pt]
    (b32.south) -- (15,-3.7) -- (8,-3.7) -- (8,0) -- (b10.west);
    \draw[box] (-0.7,0.5) rectangle (0.7,-3.5);
    \draw[box] (1.3,0.5) rectangle (2.7,-3.5);
    \draw[box] (3.3,0.5) rectangle (4.7,-3.5);
    \draw[box] (5.3,0.5) rectangle (6.7,-3.5);
    \draw[box] (8.3,0.5) rectangle (9.7,-3.5);
    \draw[box] (10.3,0.5) rectangle (11.7,-3.5);
    \draw[box] (12.3,0.5) rectangle (13.7,-3.5);
    \draw[box] (14.3,0.5) rectangle (15.7,-3.5);
    \node[lab] at (a10.south) {$\g_{Y_4,0}$};
    \node[lab] at (a20.south) {$\g_{Y_3,0}$};
    \node[lab] at (a30.south) {$\g_{Y_5,0}$};
    \node[lab] at (a11.south) {$\g_{Y_4,1}$};
    \node[lab] at (a21.south) {$\g_{Y_3,1}$};
    \node[lab] at (a31.south) {$\g_{Y_5,1}$};
    \node[lab] at (a12.south) {$\g_{Y_4,q^{m-1}-2}$};
    \node[lab] at (a22.south) {$\g_{Y_3,q^{m-1}-2}$};
    \node[lab] at (a32.south) {$\g_{Y_5,q^{m-1}-2}$};

    \node[lab] at (b10.south) {$\e_{b_i,Y_6;0}$};
    \node[lab] at (b20.south) {$\e_{b_i,Y_3;0}$};
    \node[lab] at (b30.south) {$\e_{b_i,Y_7;0}$};
    \node[lab] at (b11.south) {$\e_{b_i,Y_6;1}$};
    \node[lab] at (b21.south) {$\e_{b_i,Y_3;1}$};
    \node[lab] at (b31.south) {$\e_{b_i,Y_7;1}$};
    \node[lab] at (b12.south) {$\e_{b_i,Y_6;q^{m-1}-2}$};
    \node[lab] at (b22.south) {$\e_{b_i,Y_3;q^{m-1}-2}$};
    \node[lab] at (b32.south) {$\e_{b_i,Y_7;q^{m-1}-2}$};

    \node[tlab] at (-1.2,0) {$W$};
    \node[tlab] at (0,0.5) {$D_{Y_4}$};
    \node[tlab] at (2,0.5) {$D_{Y_3}$};
    \node[tlab] at (4,0.5) {$\cdots$};
    \node[tlab] at (6,0.5) {$D_{Y_5}$};
    \node[tlab] at (7.8,0) {$W(b_i)$};
    \node[tlab] at (9,0.5) {$E_{b_i,Y_6}$};
    \node[tlab] at (11,0.5) {$E_{b_i,Y_3}$};
    \node[tlab] at (13,0.5) {$\cdots$};
    \node[tlab] at (15,0.5) {$E_{b_i,Y_7}$};
\end{tikzpicture}
  \caption{{\footnotesize   H-path $W$ of $\G(\cup_Y D_Y)$;  and  H-path $W(b_i)$ of  $\G(\cup_Y E_{b_i,Y})$   in Lemma 4.3.{\bf (2).}}}
\vskip 2mm
  \centering
  \begin{tikzpicture}[
    scale=0.75, transform shape,
    p/.style={draw, circle, fill=white, inner sep=0pt, minimum size=2.6mm},
    d/.style={fill=white, inner sep=2pt},
    box/.style={draw, rounded corners=6pt},
    lab/.style={anchor=north, yshift=3pt, font=\footnotesize},
    tlab/.style={anchor=south, font=\small},
    midarrow/.style={
      postaction=decorate,
      decoration={markings, mark=at position 0.5 with {\arrow{angle 45}}}
    },
    remidarrow/.style={
      postaction=decorate,
      decoration={markings, mark=at position 0.5 with {\arrowreversed{angle 45}}}
    }
  ]
    \node[p] (a10)  at (0, 0)    {};
    \node[d] (d10)  at (2, 0)  {$\cdots$};
    \node[p] (a20)  at (4, 0)  {};
    \node[p] (a30)  at (6, 0)  {};
    \node[p] (c1)  at (8, 0)  {};
    \node[p] (b10)  at (10, 0)  {};
    \node[d] (d20)  at (11, 0)  {$\cdots$};
    \node[p] (b20)  at (12, 0)  {};
    \node[p] (b30)  at (14, 0)  {};
    \draw[midarrow] (a10.east) -- (d10.west);
    \draw[midarrow] (d10.east) -- (a20.west);
    \draw[midarrow] (b10.east) -- (d20.west);
    \draw[midarrow] (d20.east) -- (b20.west);
    \draw[midarrow] (b20.east) -- (b30.west);
    \node[p] (a11) at (0, -1)  {};
    \node[d] (d11) at (2, -1)  {$\cdots$};
    \node[p] (a21) at (4, -1)  {};
    \node[p] (a31) at (6, -1)  {};
    \node[p] (d3) at (7.5,-1.5) {};
    \node[d] (d) at (8, -1.5) {$\cdots$};
    \node[p] (d4) at (8.5,-1.5) {};
    \node[p] (b11)  at (10, -1.2)  {};
    \node[d] (d21)  at (11, -1.2)  {$\cdots$};
    \node[p] (b21)  at (12, -1.2)  {};
    \node[p] (b31)  at (14, -1.2)  {};
    \draw[midarrow] (a11.east) -- (d11.west);
    \draw[midarrow] (d11.east) -- (a21.west);
    \draw[midarrow] (a21.east) -- (a31.west);
    \draw[midarrow] (b11.east) -- (d21.west);
    \draw[midarrow] (d21.east) -- (b21.west);
    \draw[midarrow] (b21.east) -- (b31.west);
    \node[d] (v1) at (0, -2)   {$\vdots$};
    \node[d] (v2) at (2, -2) {$\vdots$};
    \node[d] (v3) at (4, -2) {$\vdots$};
    \node[d] (v4) at (6, -2)  {$\vdots$};
    \node[d] (v5) at (11, -2)  {$\vdots$};
    \node[d] (v6) at (14, -2)  {$\vdots$};
    \node[p] (a12) at (0, -3)   {};
    \node[d] (d12) at (2, -3) {$\cdots$};
    \node[p] (a22) at (4, -3) {};
    \node[p] (a32) at (6, -3) {};
    \node[p] (d1) at (7.5,-3) {};
    \node[d] (dd) at (8, -3) {$\cdots$};
    \node[p] (d2) at (8.5,-3) {};
    \node[p] (b12)  at (10, -3)  {};
    \node[d] (d22)  at (11, -3)  {$\cdots$};
    \node[p] (b22)  at (12, -3)  {};
    \draw[midarrow] (a12.east) -- (d12.west);
    \draw[midarrow] (d12.east) -- (a22.west);
    \draw[midarrow] (a22.east) -- (a32.west);
    \draw[midarrow] (b12.east) -- (d22.west);
    \draw[midarrow] (d22.east) -- (b22.west);

    \draw[midarrow] (a30.west)  -- (a11.east);
    \draw[midarrow] (a31.west) -- ($(a31)!0.36!(a12)$);
    \draw[midarrow] ($(a31)!0.64!(a12)$) -- (a12.east);
    \draw[midarrow, rounded corners=6pt]
      (a32.south) -- (6, -3.9) -- (-1, -3.9) -- (-1, 0) -- (a10.west);
    \draw[midarrow] (b30.west)  -- (b11.east);
    \draw[midarrow] (b31.west) -- ($(b31)!0.36!(b12)$);
    \draw[midarrow] ($(b31)!0.64!(b12)$) -- (b12.east);
    \draw[midarrow, rounded corners=6pt]
      (a20.north) -- (4, 0.7) -- (8, 0.7) -- (c1.north);
    \draw[midarrow] (c1.east)  -- (b10.west);
    \draw[midarrow, rounded corners=6pt]
      (b22.south) -- (12, -3.9) -- (8.5, -3.9) -- (d2.south);
    \draw[midarrow] (d1.east)  -- (d4.west);
    \draw[midarrow] (d3.west)  -- (a30.east);

    \draw[box] (-0.7, 0.5) rectangle (0.7, -3.7);
    \draw[box] (1.3, 0.5) rectangle (2.7, -3.7);
    \draw[box] (3.3, 0.5) rectangle (4.7, -3.7);
    \draw[box] (5.3, 0.5) rectangle (6.7, -3.7);
    \draw[box] (9.5, 0.5) rectangle (12.5, -0.5);
    \draw[box] (9.5, -0.7) rectangle (12.5, -1.7);
    \draw[box] (9.5, -2.5) rectangle (12.5, -3.5);
    \draw[box] (7, -1) rectangle (9, -2);
    \draw[box] (7, -2.5) rectangle (9, -3.5);

    \node[lab] at (a10.south)  {$\d_{Y_2,0}$};
    \node[lab] at (a20.south)  {$\d_{Y_3,0}$};
    \node[lab] at (a30.south)  {$\d_{Y_1,0}$};
    \node[lab] at (b10.south)  {$\g_{Y_3,0}$};
    \node[lab] at (b20.south)  {$\g_{Y_4,0}$};
    \node[lab] at (b30.south)  {$\b_{2,1}$};

    \node[lab] at (a11.south)  {$\d_{Y_2,1}$};
    \node[lab] at (a21.south)  {$\d_{Y_3,1}$};
    \node[lab] at (a31.south)  {$\d_{Y_1,1}$};
    \node[lab] at (b11.south)  {$\e_{b_0,Y_3;2}$};
    \node[lab] at (b21.south)  {$\e_{b_0,Y_6;2}$};
    \node[lab] at (b31.south)  {$\b_{2,3}$};

    \node[lab] at (a12.south)  {$\d_{Y_2,\frac{q-1}{2}r-1}$};
    \node[lab] at (a22.south)  {$\d_{Y_3,\frac{q-1}{2}r-1}$};
    \node[lab] at (a32.south)  {$\d_{Y_1,\frac{q-1}{2}r-1}$};
    \node[lab] at (b12.south)  {$\e_{b_{q-1},Y_3;2q}$};
    \node[lab] at (b22.south)  {$\e_{b_{q-1},Y_6;2q}$};

    \node[lab] at (c1.south)  {$\b_{2,r-1}$};
    \node[lab] at (d2.south)  {$\b_{2,2q+1}$};
    \node[lab] at (d1.south)  {$\b_{2,0}$};
    \node[lab] at (d4.south)  {$\b_{1,r-1}$};
    \node[lab] at (d3.south)  {$\b_{1,0}$};

    \node[tlab] at (0, 0.6) {$C_{Y_2}$};
    \node[tlab] at (2, 0.6) {$\cdots$};
    \node[tlab] at (4, 0.6) {$C_{Y_3}$};
    \node[tlab] at (6, 0.6) {$C_{Y_1}$};
    \node[tlab] at (10, 0.6) {$W$};
    \node[tlab] at (10, -2.3) {$W(b_0)$};
    \node[tlab] at (12, -2.6) {$W(b_{q-1})$};

    \node[tlab] at (7.5, -2.6) {$B_2'$};
    \node[tlab] at (7.5, -1.1) {$B_1$};
  \end{tikzpicture}
  \caption{{\footnotesize H-cycle $C$ of the graph $\G$ in Lemma 4.3.{\bf (2).}}}
\end{figure}

{\bf  (2)}  $q$ is even. \vskip 2mm

By Eq(\ref{minus}), set    $$Q(X)=X_1X_2^t+a_1^2+a_1a_2+\mu a_2^2\,\, {\rm and}\,\, f(X,Y)=X_1Y_2^t+X_2Y_1^t+a_1b_2+a_2b_1,$$
where  $\mu $ is given  such that  $x^2+x+\mu $ in irreducible;    $X=(X_1,X_2,a_1,a_2),$ $Y=(Y_1,Y_2,b_1,b_2)\in V(2m, q)$.
  Let $\O$ be the set of singular points so that  $|\O|=\frac{(q^m+1)(q^{m-1}-1)}{q-1}$.
     Take an $\a\in \O$. Then $\D_s=\{\b \in \O \di f(\a,\b)=0\}$, where $|\D_s|=\frac {q(q^{m-1}+1)(q^{m-2}-1)}{q-1}$.

Let $A, \rho, r, 0_{m-1}, X_0$ and $X_{e,Y}$  have the same meanings in {\bf (1)} (see Eq{(\ref{odd2})}).
Set
$$\begin{array}{llll}
B_1&=&\{\b_{1,i}\di 0\le i \le r-1\}, \quad  &\b_{1,i}=\lg(X_0A^i,0_{m-1},0,0)\rg ;\\
B_2&=&\{\b_{2,i}\di 0\le i \le r-1\}, \quad  &\b_{2,i}=\lg(0_{m-1},X_0(A^i)^\t,0,0)\rg ;\\
C_{Y}&=&\{\d_{Y, k}\di 0\le k \le q^{m-1}-2\}, &\d_{Y, k}=\lg(X_0A^k,X_{0,Y}(A^k)^\t,0,0)\rg , \\
&&& Y=\lg \th\rg Y_1, Y_1\in V(m-2,q)\setminus \{0_{m-2}\};\\
D_{Y}&=&\{\g_{Y, k}\di 0\le k \le q^{m-1}-2\}, &\g_{Y, k}=\lg(X_0A^k,X_{\mu , Y}{(A^k)^\t}, 0, 1)\rg ,\\
&&& Y\in V(m-2,q);\\
E_{c,Y}&=&\{\e_{b, Y; k}\di 0\le k \le q^{m-1}-2\}, &\e_{b, Y; k}=\lg(X_0A^k,X_{1+b+\mu b^2\th, Y}{(A^k)^\t}, 1, b)\rg ,\\
&&& Y\in V(m-2,q), b\in \FF_q.
\end{array} $$
 The remaining arguments are completely same as in {\bf (1)}, just by seting
{\small $$Y_4=Y_6(b)=Y_1=(1,0,0,\cdots,0), \, \, Y_5=Y_7(b)=Y_2=(0,1,0,\cdots,0),\, \,  Y_3=(1,1,0,\cdots,0).$$}
\qed

\section{Actions on  nonsingular points}
By Proposition~\ref{pro:3},  $M=\O(2m+1,3), \O^{\pm}(2m,3)$ or $\O^{\pm }(2m,2)$.
\subsection{ $M=\O(2m+1,3)$}
\begin{lem}\label{O41}
Let $M=\O(2m+1,3)$,  acting  an  orbit $\O$ of nongingular points.  Then  $X(\O, \Del_s)$  has a H-cycle.
\end{lem}
\demo  Since $\P\O(5,3)\cong \PSp(4,3)$, we assume $m\ge 3$.
By Eq(\ref{odd}), set  $$f(X,Y)=ab+X_1Y_2^t+X_2Y_1^t,$$ for any $X=(a,X_1,X_2)$ and $Y=(b,Y_1,Y_2)$, where $X_i,Y_i \in V(m,q)$.
Then totally  we have $3^{2m}$  nonsingular points, consisting of two $M$-orbits;
 $$\O_{\pm 1}=\{\lg X\rg \di f(X,X)=\pm 1\},\, \, {\rm where}\, \, |\O_{\pm 1}|=\frac{3^m(3^m\pm 1 )}{2}.\vspace{-2mm}$$
  For both cases,  $\D_s=\{\b\in\O_{\pm 1} \di f(\a,\b)=0\}$.
For $\O_{-1}$, we have $|\D_s|=\frac{3^{m-1}(3^m+1)}{2}\ge\frac{1}{3}|\O_{-1}|$, and so  set $\O:=\O_1$, where  $|\D_s|=\frac{3^{m-1}(3^m-1)}{2}.$

 Let $A, \rho,  r, 0_m, X_0$ and $X_{e,Y}$  have the same meanings in Lemma~\ref{os} (1) (see Eq{(\ref{odd2})}). Set
 $$\begin{array}{llll}
\{\a\},  &&\a=\lg(1,0_m,0_m)\rg;\\
B&=&\{\b_k \di 0\le k \le 3^m-2\},\quad &\b_k=\lg(1,0_m,X_0(A^k)^\t)\rg;\\
C_Y&=&\{\d_{Y,k} \di 0 \le k \le \frac{3^m-1}{2}-1\},\quad & \d_{Y,k}=\lg(0,X_0A^k,X_{2,Y}(A^k)^\t)\rg,\\
&&& Y\in V(m-1,3);\\
D_Y&=&\{\g_{Y,k} \di 0 \le k \le 3^m-2\},\quad &\g_{Y,k}=\lg(1,X_0A^k,X_{0,Y}(A^k)^\t)\rg,\\
&&& Y\in V(m-1,3).
\end{array}$$
  Check that $\O$ is a disjoint union of the above  $\lg\rho\rg$-orbits, that is
  $$\O=\{\a\}\cup B\cup (\cup_{Y}C_Y)\cup (\cup_{Y}D_Y),$$
where $Y\in V(m-1,3).$.   Label    all  vectors  of $V(m-3,3)$ by  $\{ Y_j\di 0\le j \le 3^{m-3}-1\}$ and   write
  every   vector  $Y\in V(m-1, 3)$  by $Y=(y_1,y_2, Y_j)$, where $y_1, y_2\in \FF_q$. Then we carry the proof by the following three steps.
\vskip 2mm

{\bf Step 1:}  Check that  $\g_{(y_1,y_2, Y_i), 0}\thicksim \g_{(y_1',y_2', Y_j), 0}^{\rho }$ if and only if   $1+y_1'+y_2=0$, where $0\le i,j \le 3^{m-3}-1$.
Given any $Y_j\in V(m-3,3)$,  there is  a path:
$$\begin{array}{ll} U_0(Y_j):\,  &\g_{(0, 0, Y_j), 0} \thicksim \g_{(2, 0, Y_j), 1} \thicksim \g_{(2, 1, Y_j), 2}\thicksim \g_{(1, 0, Y_j), 3} \thicksim\g_{(2, 2, Y_j), 4} \thicksim \g_{(0, 1, Y_j), 5} \\
&\thicksim \g_{(1, 1, Y_j), 6}  \thicksim \g_{(1, 2, Y_j), 7}\thicksim \g_{(0, 2, Y_j), 8}.\end{array} $$
Since  $\lg \rho ^9\rg =\lg \rho \rg $,  we get a path (see the left part of Figure 4):
 $$U(Y_j):\quad U_0(Y_j), U_0(Y_j)^{\rho ^9},  U_0(Y_j)^{\rho ^{9\cdot 2}}, \cdots ,   U_0(Y_j)^{\rho ^{9(3^m-2)}}.$$
Then  we get a H-cycle $U$ for $\G(\cup_Y D_Y)$:
$$U:\quad U(Y_0),U(Y_1),U(Y_2),\cdots,U(Y_{3^{m-3}-1}),\g_{(0,0,Y_0),0}.$$
Moreover,  set $9k=-2(\mod 3^m-1)$. Then   $U_0(Y_0)^{\rho ^{9\cdot k}}$ passes  the edge   $\{\g_{(2,1,Y_0),0}, $ $\g_{(1,0,Y_0),1}\}$.
Therefore, from the  H-cycle $U$, we have H-path $\g_{(2,1,Y_0),0}\to \g_{(1,0,Y_0),1}$  of $\G(\cup_Y D_Y)$.
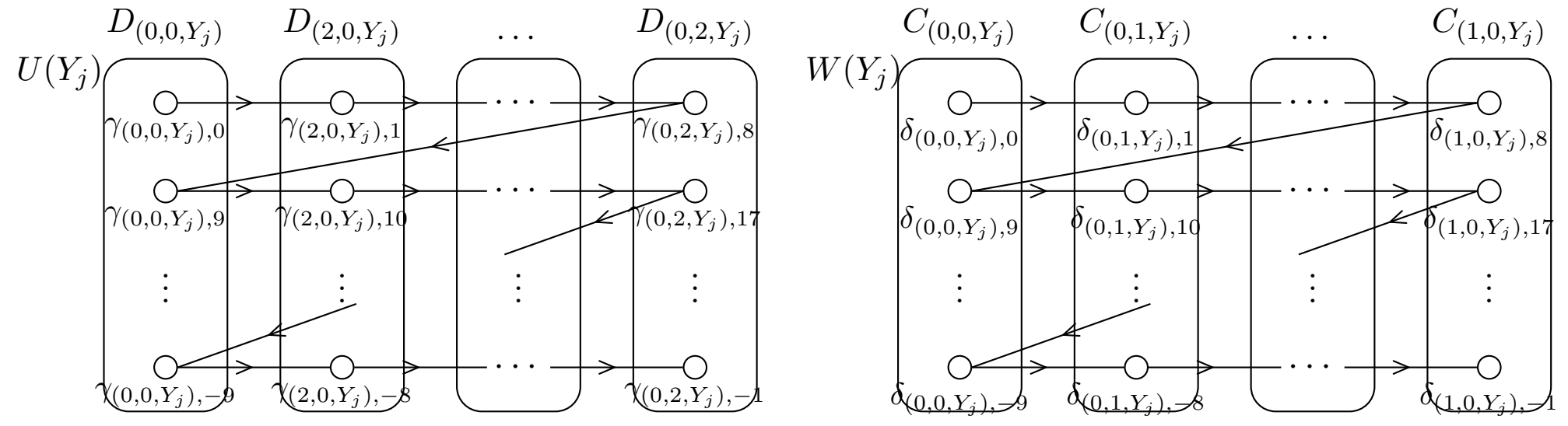
\begin{figure}[htbp]
  \centering
  \begin{tikzpicture}
  [
    scale=0.75, transform shape,
    p/.style={draw, circle, fill=white, inner sep=0pt, minimum size=2.6mm},
    d/.style={fill=white, inner sep=2pt},
    box/.style={draw, rounded corners=6pt},
    lab/.style={anchor=north, yshift=3pt, font=\footnotesize},
    tlab/.style={anchor=south, font=\small},
    midarrow/.style={
      postaction=decorate,
      decoration={markings, mark=at position 0.5 with {\arrow{angle 45}}}
    },
    remidarrow/.style={
      postaction=decorate,
      decoration={markings, mark=at position 0.5 with {\arrowreversed{angle 45}}}
    }
  ]
    \node[p] (a10)  at (0, 0)  {};
    \node[p] (a20)  at (2, 0)  {};
    \node[d] (d10)  at (4, 0)  {$\cdots$};
    \node[p] (a30)  at (6, 0)  {};
    \node[p] (b10)  at (9, 0)  {};
    \node[p] (b20)  at (11, 0)  {};
    \node[d] (d20)  at (13, 0)  {$\cdots$};
    \node[p] (b30)  at (15, 0)  {};
    \draw[midarrow] (a10.east) -- (a20.west);
    \draw[midarrow] (a20.east) -- (d10.west);
    \draw[midarrow] (d10.east) -- (a30.west);
    \draw[midarrow] (b10.east) -- (b20.west);
    \draw[midarrow] (b20.east) -- (d20.west);
    \draw[midarrow] (d20.east) -- (b30.west);
    \node[p] (a11)  at (0, -1)  {};
    \node[p] (a21)  at (2, -1)  {};
    \node[d] (d11)  at (4, -1)  {$\cdots$};
    \node[p] (a31)  at (6, -1)  {};
    \node[p] (b11)  at (9, -1)  {};
    \node[p] (b21)  at (11, -1)  {};
    \node[d] (d21)  at (13, -1)  {$\cdots$};
    \node[p] (b31)  at (15, -1)  {};
    \draw[midarrow] (a11.east) -- (a21.west);
    \draw[midarrow] (a21.east) -- (d11.west);
    \draw[midarrow] (d11.east) -- (a31.west);
    \draw[midarrow] (b11.east) -- (b21.west);
    \draw[midarrow] (b21.east) -- (d21.west);
    \draw[midarrow] (d21.east) -- (b31.west);
    \node[d] (v1)  at (0, -2)  {$\vdots$};
    \node[d] (v2)  at (2, -2)  {$\vdots$};
    \node[d] (v3)  at (4, -2)  {$\vdots$};
    \node[d] (v4)  at (6, -2)  {$\vdots$};
    \node[d] (v5)  at (9, -2) {$\vdots$};
    \node[d] (v6)  at (11, -2) {$\vdots$};
    \node[d] (v7)  at (13, -2) {$\vdots$};
    \node[d] (v8)  at (15, -2) {$\vdots$};
    \node[p] (a12)  at (0, -3)  {};
    \node[p] (a22)  at (2, -3)  {};
    \node[d] (d12)  at (4, -3)  {$\cdots$};
    \node[p] (a32)  at (6, -3)  {};
    \node[p] (b12)  at (9, -3)  {};
    \node[p] (b22)  at (11, -3)  {};
    \node[d] (d22)  at (13, -3)  {$\cdots$};
    \node[p] (b32)  at (15, -3)  {};
    \draw[midarrow] (a12.east) -- (a22.west);
    \draw[midarrow] (a22.east) -- (d12.west);
    \draw[midarrow] (d12.east) -- (a32.west);
    \draw[midarrow] (b12.east) -- (b22.west);
    \draw[midarrow] (b22.east) -- (d22.west);
    \draw[midarrow] (d22.east) -- (b32.west);

    \draw[midarrow] (a30.west) -- (a11.east);
    \draw[midarrow] (a31.west) -- ($(a31)!0.36!(a12)$);
    \draw[midarrow] ($(a31)!0.64!(a12)$) -- (a12.east);
    \draw[midarrow] (b30.west) -- (b11.east);
    \draw[midarrow] (b31.west) -- ($(b31)!0.36!(b12)$);
    \draw[midarrow] ($(b31)!0.64!(b12)$) -- (b12.east);
    \draw[box] (-0.7,0.5) rectangle (0.7,-3.5);
    \draw[box] (1.3,0.5) rectangle (2.7,-3.5);
    \draw[box] (3.3,0.5) rectangle (4.7,-3.5);
    \draw[box] (5.3,0.5) rectangle (6.7,-3.5);
    \draw[box] (8.3,0.5) rectangle (9.7,-3.5);
    \draw[box] (10.3,0.5) rectangle (11.7,-3.5);
    \draw[box] (12.3,0.5) rectangle (13.7,-3.5);
    \draw[box] (14.3,0.5) rectangle (15.7,-3.5);
    \node[lab] at (a10.south) {$\g_{(0,0,Y_j),0}$};
    \node[lab] at (a20.south) {$\g_{(2,0,Y_j),1}$};
    \node[lab] at (a30.south) {$\g_{(0,2,Y_j),8}$};
    \node[lab] at (a11.south) {$\g_{(0,0,Y_j),9}$};
    \node[lab] at (a21.south) {$\g_{(2,0,Y_j),10}$};
    \node[lab] at (a31.south) {$\g_{(0,2,Y_j),17}$};
    \node[lab] at (a12.south) {$\g_{(0,0,Y_j),-9}$};
    \node[lab] at (a22.south) {$\g_{(2,0,Y_j),-8}$};
    \node[lab] at (a32.south) {$\g_{(0,2,Y_j),-1}$};

    \node[lab] at (b10.south) {$\d_{(0,0,Y_j),0}$};
    \node[lab] at (b20.south) {$\d_{(0,1,Y_j),1}$};
    \node[lab] at (b30.south) {$\d_{(1,0,Y_j),8}$};
    \node[lab] at (b11.south) {$\d_{(0,0,Y_j),9}$};
    \node[lab] at (b21.south) {$\d_{(0,1,Y_j),10}$};
    \node[lab] at (b31.south) {$\d_{(1,0,Y_j),17}$};
    \node[lab] at (b12.south) {$\d_{(0,0,Y_j),-9}$};
    \node[lab] at (b22.south) {$\d_{(0,1,Y_j),-8}$};
    \node[lab] at (b32.south) {$\d_{(1,0,Y_j),-1}$};

    \node[tlab] at (-1.2,0) {$U(Y_j)$};
    \node[tlab] at (0,0.5) {$D_{(0,0,Y_j)}$};
    \node[tlab] at (2,0.5) {$D_{(2,0,Y_j)}$};
    \node[tlab] at (4,0.5) {$\cdots$};
    \node[tlab] at (6,0.5) {$D_{(0,2,Y_j)}$};
    \node[tlab] at (7.8,0) {$W(Y_j)$};
    \node[tlab] at (9,0.5) {$C_{(0,0,Y_j)}$};
    \node[tlab] at (11,0.5) {$C_{(0,1,Y_j)}$};
    \node[tlab] at (13,0.5) {$\cdots$};
    \node[tlab] at (15,0.5) {$C_{(1,0,Y_j)}$};
  \end{tikzpicture}
  \caption{{\footnotesize Paths $U(Y_j)$ and $W(Y_j)$ in the graph $\G$ in Lemma 5.1.}}
\end{figure}	
\vskip 2mm
{\bf Step 2:}  Check that $\d_{(y_1,y_2,Y_i),0}\thicksim \d_{(y_1',y_2',Y_j),1}$ if and only if $y_1'+y_2=0$, where $0\le i,j \le 3^{m-3}-1$. Given a $Y_j\in V(m-3,3)$, there is a path:
\begin{equation} \label{W0}\begin{array}{ll} W_0(Y_j):\,&\d_{(0,0,Y_j),0}\thicksim \d_{(0,1,Y_j),1}\thicksim \d_{(2,0,Y_j),2}\thicksim\d_{(0,2,Y_j),3} \thicksim\d_{(1,1,Y_j),4}\thicksim\d_{(2,1,Y_j),5}\\
&\thicksim\d_{(2,2,Y_j),6}\thicksim\d_{(1,2,Y_j),7}\thicksim\d_{(1,0,Y_j),8}.
\end{array}\end{equation}
Since $\lg \o \rho^9\rg=\lg\o\rho\rg$, we get a path (see the right of Figure 4):
\begin{equation}\label{W1}    W(Y_j):\quad W_0(Y_j),W_0(Y_j)^{\orho^9},W_0(Y_j)^{\orho^{9\cdot 2}},\cdots,W_0(Y_j)^{\orho^{9(\frac{3^m-1}{2}-1)}}.\end{equation}
Note that $\b_k$ is adjacent to $\d_{Y,k\pm1}$ for any given $k$. For $W(Y_0)$, since $|B|=2|D_{Y}|$, we have a new path:
$$\begin{array}{ll} R_0:\, &\d_{(0,0,Y_0),0}\thicksim \b_1 \thicksim \d_{(0,1,Y_0),2}\thicksim \b_{1+\frac{3^m-1}{2}}\thicksim\d_{(2,0,Y_0),0}\thicksim\d_{(0,2,Y_0),1} \thicksim\d_{(1,1,Y_0),2}\\
&\thicksim\d_{(2,1,Y_0),3}\thicksim\d_{(2,2,Y_0),4}\thicksim\d_{(1,2,Y_0),5}\thicksim\d_{(1,0,Y_0),6}.
\end{array}$$
We need to deal with two cases, separately.
\vskip 2mm
(i) Suppose $7 \nmid (3^m-1)$. Noting  $\lg\orho^7\rg=\lg\orho\rg$,  we have  a path (see Figure 5):
$$R:\quad R_0,R_0^{\orho^7},R_0^{\orho^{7\cdot 2}},\cdots,  R_0^{\orho^{7(\frac{3^m-1}{2}-1)}}.$$
Therefore,  we get a H-cycle $W$ of  $\G(\cup_Y C_Y\cup B)$:
$$W:\quad R,W(Y_1),W(Y_2),\cdots,W(Y_{3^{m-3}-1}),\d_{(0,0,Y_0),0}.$$

(ii) Suppose   $7 \di (3^m-1)$. Then  $6\di m$. Consider the cycle (see Figure 5):
$$S:\quad R_0,R_0^{\orho^7},R_0^{\orho^{7\cdot 2}}, \cdots , R_0^{\orho^{7(\frac{3^m-1}{14}-1)}},\d_{(0,0,Y_0),0},$$
and take  a H-path $\d_{(0,0,Y_0),0}\to \b_1$ in  $S$.
Correspondingly, we have other six cycles
$S^{\orho^j}$ and H-paths $\d_{(0,0,Y_0),0}^{\orho^j}\to\b_1^{\orho^j}$ in $S^{\orho^j}$, where $1\le j\le 6$.
Remind that $W(Y_0)\cup B$ is a disjoint union of these seven cycles $S^{\orho^j}$.

Let $W(Y_1)$ be defined in Eq(\ref{W1}).   Now we are constructing a H-path, say $W_1$, of  $\G(W(Y_1)\cup W(Y_0)\cup B)$:
first insert the H-path  $\d_{(0,0,Y_0),0}\to \b_1$  of $S$ in    $W_0(Y_1)$  (see \ref{W0}) and get a path (see Figure 6):
$$\begin{array}{ll} S(Y_1):\, &\d_{(0,0,Y_1),0}\thicksim\b_1\to\d_{(0,0,Y_0),0}\thicksim \d_{(0,1,Y_1),1}\thicksim \d_{(2,0,Y_1),2}\thicksim\d_{(0,2,Y_1),3}\\
&\thicksim\d_{(1,1,Y_1),4}\thicksim\d_{(2,1,Y_1),5}\thicksim\d_{(2,2,Y_1),6}\thicksim\d_{(1,2,Y_1),7}\thicksim\d_{(1,0,Y_1),8}.
\end{array}$$
Set $9k_j=j (\mod \frac{3^m-1}{2})$, where $0\le j\le 6$.  Then we get seven  paths $S(Y_1)^{\orho^{9\cdot k_j}}$.  Then we  get a H-path $W_1$ of  $\G(W(Y_1)\cup W(Y_0)\cup B)$:
$$W_1:\quad S(Y_1),W_0(Y_1)^{\orho^9},\cdots,S(Y_1)^{\orho^{9\cdot k_1}},\cdots,S(Y_1)^{\orho^{9\cdot k_2}},\cdots,W(Y _1)^{\orho^{9(\frac{3^m-1}{2}-1)}}.$$
Therefore, we get a H-cycle $W$ for $\G(\cup_Y C_Y\cup B)$:
$$W:\quad W_1, W(Y_2),\cdots,W(Y_{3^{m-3}-1}),\d_{(0,0,Y_0),0}.\vspace{-2mm}$$
\begin{figure}[htbp]
  \centering
  \begin{tikzpicture}
  [
    scale=0.75, transform shape,
    p/.style={draw, circle, fill=white, inner sep=0pt, minimum size=2.6mm},
    d/.style={fill=white, inner sep=2pt},
    box/.style={draw, rounded corners=6pt},
    lab/.style={anchor=north, yshift=3pt, font=\footnotesize},
    tlab/.style={anchor=south, font=\small},
    midarrow/.style={
      postaction=decorate,
      decoration={markings, mark=at position 0.5 with {\arrow{angle 45}}}
    },
    remidarrow/.style={
      postaction=decorate,
      decoration={markings, mark=at position 0.5 with {\arrowreversed{angle 45}}}
    }
  ]
    \node[p] (a10)  at (0, 0)  {};
    \node[p] (a20)  at (2, 0)  {};
    \node[p] (a30)  at (4, 0)  {};
    \node[p] (a40)  at (6, 0)  {};
    \node[p] (a50)  at (8, 0)  {};
    \node[d] (d10)  at (10, 0)  {$\cdots$};
    \node[p] (a60)  at (12, 0)  {};
    \draw[midarrow] (a10.east) -- (a20.west);
    \draw[midarrow] (a20.east) -- (a30.west);
    \draw[midarrow] (a30.east) -- (a40.west);
    \draw[midarrow] (a40.east) -- (a50.west);
    \draw[midarrow] (a50.east) -- (d10.west);
    \draw[midarrow] (d10.east) -- (a60.west);

    \node[p] (a11)  at (0, -1)  {};
    \node[p] (a21)  at (2, -1)  {};
    \node[p] (a31)  at (4, -1)  {};
    \node[p] (a41)  at (6, -1)  {};
    \node[p] (a51)  at (8, -1)  {};
    \node[d] (d11)  at (10, -1)  {$\cdots$};
    \node[p] (a61)  at (12, -1)  {};
    \draw[midarrow] (a11.east) -- (a21.west);
    \draw[midarrow] (a21.east) -- (a31.west);
    \draw[midarrow] (a31.east) -- (a41.west);
    \draw[midarrow] (a41.east) -- (a51.west);
    \draw[midarrow] (a51.east) -- (d11.west);
    \draw[midarrow] (d11.east) -- (a61.west);
    \node[d] (v1)  at (0, -2)  {$\vdots$};
    \node[d] (v2)  at (2, -2)  {$\vdots$};
    \node[d] (v3)  at (4, -2)  {$\vdots$};
    \node[d] (v4)  at (6, -2)  {$\vdots$};
    \node[d] (v5)  at (8, -2)  {$\vdots$};
    \node[d] (v6)  at (10, -2)  {$\vdots$};
    \node[d] (v7)  at (12, -2)  {$\vdots$};
    \node[p] (a12)  at (0, -3)  {};
    \node[p] (a22)  at (2, -3)  {};
    \node[p] (a32)  at (4, -3)  {};
    \node[p] (a42)  at (6, -3)  {};
    \node[p] (a52)  at (8, -3)  {};
    \node[d] (d12)  at (10, -3)  {$\cdots$};
    \node[p] (a62)  at (12, -3)  {};
    \draw[midarrow] (a12.east) -- (a22.west);
    \draw[midarrow] (a22.east) -- (a32.west);
    \draw[midarrow] (a32.east) -- (a42.west);
    \draw[midarrow] (a42.east) -- (a52.west);
    \draw[midarrow] (a52.east) -- (d12.west);
    \draw[midarrow] (d12.east) -- (a62.west);

    \draw[midarrow] (a60.west) -- (a11.east);
    \draw[midarrow] (a61.west) -- ($(a61)!0.36!(a12)$);
    \draw[midarrow] ($(a61)!0.64!(a12)$) -- (a12.east);
    \draw[midarrow, dashed,  rounded corners=6pt]
    (a62.south) -- (12,-3.8) -- (-0.9,-3.8) -- (-0.9, 0) -- (a10.west);
    \draw[box] (-0.7,0.5) rectangle (0.7,-3.6);
    \draw[box] (1.3,0.5) rectangle (2.7,-3.6);
    \draw[box] (3.3,0.5) rectangle (4.7,-3.6);
    \draw[box] (5.3,0.5) rectangle (6.7,-3.6);
    \draw[box] (7.3,0.5) rectangle (8.7,-3.6);
    \draw[box] (9.3,0.5) rectangle (10.7,-3.6);
    \draw[box] (11.3,0.5) rectangle (12.7,-3.6);

    \node[lab] at (a10.south) {$\d_{(0,0,Y_0),0}$};
    \node[lab] at (a20.south) {$\b_1$};
    \node[lab] at (a30.south) {$\d_{(0,1,Y_0),2}$};
    \node[lab] at (a40.south) {$\b_{1+\frac{3^m-1}{2}}$};
    \node[lab] at (a50.south) {$\d_{(2,0,Y_0),0}$};
    \node[lab] at (a60.south) {$\d_{(1,0,Y_0),6}$};
    \node[lab] at (a11.south) {$\d_{(0,0,Y_0),7}$};
    \node[lab] at (a21.south) {$\b_8$};
    \node[lab] at (a31.south) {$\d_{(0,1,Y_0),9}$};
    \node[lab] at (a41.south) {$\b_{8+\frac{3^m-1}{2}}$};
    \node[lab] at (a51.south) {$\d_{(2,0,Y_0),7}$};
    \node[lab] at (a61.south) {$\d_{(1,0,Y_0),13}$};

    \node[lab] at (a12.south) {$\d_{(0,0,Y_0),-7}$};
    \node[lab] at (a22.south) {$\b_{-6}$};
    \node[lab] at (a32.south) {$\d_{(0,1,Y_0),-5}$};
    \node[lab] at (a42.south) {$\b_{-6+\frac{3^m-1}{2}}$};
    \node[lab] at (a52.south) {$\d_{(2,0,Y_0),-7}$};
    \node[lab] at (a62.south) {$\d_{(1,0,Y_0),-1}$};

    \node[tlab] at (-1.2,0.5) {$R/S$};
    \node[tlab] at (-1.2,-0.3) {$R_0$};
    \node[tlab] at (-1.2,-1.3) {$R_0^{\orho^7}$};
    \node[tlab] at (-1.2,-2.3) {$\vdots$};
    \node[tlab] at (-1.3,-3.3) {$R_0^{\orho^{-7}}$};
    \node[tlab] at (0,0.5) {$C_{(0,0,Y_0)}$};
    \node[tlab] at (4,0.5) {$C_{(0,1,Y_0)}$};
    \node[tlab] at (8,0.5) {$C_{(2,0,Y_0)}$};
    \node[tlab] at (10,0.5) {$\cdots$};
    \node[tlab] at (12,0.5) {$C_{(1,0,Y_0)}$};
  \end{tikzpicture}
  \caption{{\footnotesize Path $R$ or Cycle $S$ in $\G$ in Lemma 5.1}}

\vskip 1.5mm
  \begin{tikzpicture}
  [
    scale=0.75, transform shape,
    p/.style={draw, circle, fill=white, inner sep=0pt, minimum size=2.6mm},
    d/.style={fill=white, inner sep=2pt},
    box/.style={draw, rounded corners=6pt},
    lab/.style={anchor=north, yshift=3pt, font=\footnotesize},
    tlab/.style={anchor=south, font=\small},
    midarrow/.style={
      postaction=decorate,
      decoration={markings, mark=at position 0.5 with {\arrow{angle 45}}}
    },
    remidarrow/.style={
      postaction=decorate,
      decoration={markings, mark=at position 0.5 with {\arrowreversed{angle 45}}}
    }
  ]
    \node[p] (a10)  at (0, 0)  {};
    \node[p] (a20)  at (2, 0)  {};
    \node[p] (a30)  at (4, 0)  {};
    \node[p] (a40)  at (6, 0)  {};
    \node[p] (a50)  at (8, 0)  {};
    \node[d] (d10)  at (10, 0)  {$\cdots$};
    \node[p] (a60)  at (12, 0)  {};
    \draw[midarrow] (a10.east) -- (a20.west);
    \draw[midarrow] (a20.east) -- (a30.west);
    \draw[midarrow] (a30.east) -- (a40.west);
    \draw[midarrow] (a40.east) -- (a50.west);
    \draw[midarrow] (a50.east) -- (d10.west);
    \draw[midarrow] (d10.east) -- (a60.west);

    \node[p] (a11)  at (0, -1)  {};
    \node[p] (a41)  at (6, -1)  {};
    \node[p] (a51)  at (8, -1)  {};
    \node[d] (d11)  at (10, -1)  {$\cdots$};
    \node[p] (a61)  at (12, -1)  {};
    \draw[midarrow] (a11.east) -- (a41.west);
    \draw[midarrow] (a41.east) -- (a51.west);
    \draw[midarrow] (a51.east) -- (d11.west);
    \draw[midarrow] (d11.east) -- (a61.west);
    \node[d] (v1)  at (0, -2)  {$\vdots$};
    \node[d] (v2)  at (3, -2)  {$\vdots$};
    \node[d] (v4)  at (6, -2)  {$\vdots$};
    \node[d] (v5)  at (8, -2)  {$\vdots$};
    \node[d] (v6)  at (10, -2)  {$\vdots$};
    \node[d] (v7)  at (12, -2)  {$\vdots$};
    \node[p] (a12)  at (0, -3)  {};
    \node[p] (a22)  at (2, -3)  {};
    \node[p] (a32)  at (4, -3)  {};
    \node[p] (a42)  at (6, -3)  {};
    \node[p] (a52)  at (8, -3)  {};
    \node[d] (d12)  at (10, -3)  {$\cdots$};
    \node[p] (a62)  at (12, -3)  {};
    \draw[midarrow] (a12.east) -- (a22.west);
    \draw[midarrow] (a22.east) -- (a32.west);
    \draw[midarrow] (a32.east) -- (a42.west);
    \draw[midarrow] (a42.east) -- (a52.west);
    \draw[midarrow] (a52.east) -- (d12.west);
    \draw[midarrow] (d12.east) -- (a62.west);
    \node[d] (v1)  at (0, -4)  {$\vdots$};
    \node[d] (v2)  at (3, -4)  {$\vdots$};
    \node[d] (v4)  at (6, -4)  {$\vdots$};
    \node[d] (v5)  at (8, -4)  {$\vdots$};
    \node[d] (v6)  at (10, -4)  {$\vdots$};
    \node[d] (v7)  at (12, -4)  {$\vdots$};
    \node[p] (a13)  at (0, -5)  {};
    \node[p] (a43)  at (6, -5)  {};
    \node[p] (a53)  at (8, -5)  {};
    \node[d] (d13)  at (10, -5)  {$\cdots$};
    \node[p] (a63)  at (12, -5)  {};
    \draw[midarrow] (a13.east) -- (a43.west);
    \draw[midarrow] (a43.east) -- (a53.west);
    \draw[midarrow] (a53.east) -- (d13.west);
    \draw[midarrow] (d13.east) -- (a63.west);

    \draw[midarrow] (a60.west) -- (a11.east);
    \draw[midarrow] (a61.west) -- ($(a61)!0.36!(a12)$);
    \draw[midarrow] ($(a61)!0.64!(a12)$) -- (a12.east);
    \draw[midarrow] (a62.west) -- ($(a62)!0.36!(a13)$);
    \draw[midarrow] ($(a62)!0.64!(a13)$) -- (a13.east);
    \draw[box] (-0.7,0.5) rectangle (0.7,-5.6);
    \draw[box] (5.3,0.5) rectangle (6.7,-5.6);
    \draw[box] (7.3,0.5) rectangle (8.7,-5.6);
    \draw[box] (9.3,0.5) rectangle (10.7,-5.6);
    \draw[box] (11.3,0.5) rectangle (12.7,-5.6);
    \draw[box] (1.5,0.5) rectangle (4.5,-0.6);
    \draw[box] (1.5,-2.5) rectangle (4.5,-3.6);


    \node[lab] at (a10.south) {$\d_{(0,0,Y_1),0}$};
    \node[lab] at (a20.south) {$\b_1$};
    \node[lab] at (a30.south) {$\d_{(0,0,Y_0),0}$};
    \node[lab] at (a40.south) {$\d_{(0,1,Y_1),1}$};
    \node[lab] at (a50.south) {$\d_{(2,0,Y_1),2}$};
    \node[lab] at (a60.south) {$\d_{(1,0,Y_1),8}$};
    \node[lab] at (a11.south) {$\d_{(0,0,Y_1),9}$};
    \node[lab] at (a41.south) {$\d_{(0,1,Y_1),10}$};
    \node[lab] at (a51.south) {$\d_{(2,0,Y_1),11}$};
    \node[lab] at (a61.south) {$\d_{(1,0,Y_1),17}$};
    \node[lab] at (a12.south) {$\d_{(0,0,Y_1),1}$};
    \node[lab] at (a22.south) {$\b_2$};
    \node[lab] at (a32.south) {$\d_{(0,0,Y_0),1}$};
    \node[lab] at (a42.south) {$\d_{(0,1,Y_1),2}$};
    \node[lab] at (a52.south) {$\d_{(2,0,Y_1),3}$};
    \node[lab] at (a62.south) {$\d_{(1,0,Y_1),9}$};
    \node[lab] at (a13.south) {$\d_{(0,0,Y_1),-9}$};
    \node[lab] at (a43.south) {$\d_{(0,1,Y_1),-8}$};
    \node[lab] at (a53.south) {$\d_{(2,0,Y_1),-7}$};
    \node[lab] at (a63.south) {$\d_{(1,0,Y_1),-1}$};

    \node[tlab] at (-1.2,0.5) {$W_1$};
    \node[tlab] at (-1.2,-0.3) {$S(Y_1)$};
    \node[tlab] at (-1.5,-1.3) {$W_0(Y_1)^{\orho^9}$};
    \node[tlab] at (-1.2,-2.3) {$\vdots$};
    \node[tlab] at (-1.3,-3.3) {$S(Y_1)^{\orho}$};
    \node[tlab] at (-1.6,-5.3) {$W_0(Y_1)^{\orho^{-9}}$};
    \node[tlab] at (0,0.5) {$C_{(0,0,Y_1)}$};
    \node[tlab] at (6,0.5) {$C_{(0,1,Y_1)}$};
    \node[tlab] at (8,0.5) {$C_{(2,0,Y_1)}$};
    \node[tlab] at (10,0.5) {$\cdots$};
    \node[tlab] at (12,0.5) {$C_{(0,2,Y_1)}$};
    \node[tlab] at (2,0.5) {$S$};
    \node[tlab] at (2,-2.5) {$S^{\orho}$};
  \end{tikzpicture}
  \caption{{\footnotesize Paths $S(Y_1)$ and $W_1$ in $\G$ in Lemma 5.1}}
\end{figure}

\vskip 2mm
{\bf Step 3:}  Pick a H-path   $\d_{(2,0,Y_0),1}\to \d_{(0,2,Y_0),2}$ in cycle $W$, for both cases $7\di (3^m-1)$ and $7\nmid (3^m-1)$.  Since  $\a$ is adjacent to all points in $C_Y$ for all $Y$, we get a  H-path
$\d_{(2,0,Y_0),1} \to \a,\d_{(0,2,Y_0),2}$ of  $\G(\cup_Y C_Y\cup B\cup \{\a\})$.
By attaching the H-path of $\G(\cup_Y D_Y)$ in {\bf Step 1},     we  finally can get a H-cycle $C$ for $\G$:
$$C:\quad \d_{(2,0,Y_0),1} \to \a,\d_{(0,2,Y_0),2},   \g_{(1,0,Y_0),1}\to \g_{(2,1,Y_0),0}, \d_{(2,0,Y_0);1}.\vspace{-2mm}$$
\subsection{ $M=\O^{\pm}(2m,3)$ and $\O^{\pm }(2m,2)$}
\begin{lem}\label{O42}
Let $M=\O^{\pm}(2m,3)$ and $\O^{\pm }(2m,2)$,   acting  an  orbit $\O$ of nongingular points.  Then  $X(\O, \Del_s)$  has a H-cycle.
\end{lem}
\demo First we let  $M=\O^{+}(2m,3).$
By Eq(\ref{plus}), set    $$f(X,Y)=X_1Y_2^t+X_2Y_1^t,  \, \forall X=(X_1,X_2),Y=(Y_1,Y_2),  X_i,Y_i \in V(m,q).$$
 Then totally  we have  $3^{m-1}(3^m-1)$  nonsingular points, consisting of two $M$-orbits,
   $\O_{\pm 1}=\{\lg X\rg \di f(X,X)=\pm 1\}$,   where $|\O_{\pm 1}|=\frac{3^{m-1}(3^m-1)}{2}$.
It suffices to consider $\O=\O_{-1}$.  For any $\a\in\O$, $\D_s=\{\b\in\O|f(\a,\b)=0\}$, where $|\D_s|=\frac{3^{m-1}(3^{m-1}-1)}{2}$.

Let $A, \rho,  r, 0_m, X_0$ and $X_{e,Y}$  have the same meanings in Lemma~\ref{esp} (1) (see Eq{(\ref{odd2})}).
Check that all $\lg\rho\rg$-orbits on $\O$ are
$C_Y=\{\d_{Y,k}\di 0\le k \le \frac{3^m-1}{2}-1\}$, where
$$\d_{Y,k}=\lg(X_0A^k,X_{1,Y}(A^k)^\t)\rg,  \quad Y_1\in V(m-1,q).$$
Check that  $\d_{(y_1,y_2, Y_i), 0}\thicksim \d_{(y_1',y_2', Y_j), 0}^{\rho }$ if and only if   $y_1'+y_2=0$, where $0\le i,j\le 3^{m-3}-1$.
Given any $Y_j\in V(m-3,3)$,  we get a path:
$$\begin{array}{ll} W_0(Y_j):\,&\d_{(0, 0, Y_j), 0}\thicksim \d_{(0, 1, Y_j), 1} \thicksim \d_{(2, 0, Y_j), 2} \thicksim
 \d_{(0, 2, Y_j), 3} \thicksim \d_{(1, 1, Y_j), 4}\thicksim\d_{(2, 1, Y_j), 5}\\
&\thicksim \d_{(2, 2, Y_j), 6} \thicksim \d_{(1, 2, Y_j), 7} \thicksim\d_{(1, 0, Y_j), 8}.\end{array}$$
Since $\lg\orho^9\rg=\lg\orho\rg$, we get a path:
$$W(Y_j):\quad W_0(Y_j),W_0(Y_j)^{\orho^9},W_0(Y_j)^{\orho^{9\cdot 2}},\cdots,W_0(Y_j)^{\orho^{9(\frac{3^m-1}{2}-1)}}$$
for any $j$ and then  get a H-cycle $C$ of $\G$:
$$C:\quad W(Y_0),W(Y_1),W(Y_2)\cdots,W(Y_{3^{m-3}-1}),\d_{(0,0,Y_0),0}.$$

Second, let  $M=\O^{\pm }(2m,2)$ and $\O^{-}(2m,3)$. Then
 one may  easily check  that $|\D_s|\ge |\O|/3$, for all three cases.
\qed

\section{Actions on  maximal totally  singular  subspaces}
 By  Proposition~\ref{pro:3}, $M\in \{ \O^-(6,q), \O^{+}(8,q), \O^{+}(10,q)\}$.  Since   $\O^-(6,q)\cong \SU(4,q)$, which has been discussed in \cite{DYao26}, we  just deal with the other  two groups.

A  plane  $\lg u, v\rg$  of $V$  is said to be {\it hyperbolic } if $f(u,v)=1$ and $Q(u)=Q(v)=0$ so that $f(u,u)=f(v,v)=0$, no matter $q$ is odd or even.
We shall use abbreviations  t.s. and m.t.s. respectively for  totally singular  and   maximal totally singular.

\subsection{$M=\O^{+}(8,q)$}
Let $V=\lg u_1,u_2\rg\perp \lg v_1,v_2\rg\perp \lg w_1,w_2\rg \perp \lg x_1,x_2\rg$ be a  decomposition  of   hyperbolic planes for $V$.
Let    $U=\lg u_1,v_1, w_1, x_1\rg$  be an m.t.s. subspace of $V$.
The group $M:=\GO^+(8,q)$ is transitive on the set of m.t.s. subspaces of $V$, consisting of    two $M$-orbits,  say $\O_1$ and $\O_2$, while $\O_1$ (resp. $\O_2$)   are
   all those  intersecting $U$ with even  (resp. odd) dimension. Since the actions of $M$ on $\O_1$ and $\O_2$ are equivalent,
we just  consider its action on  $\O=\O_1$, where $|\O|=(q+1)(q^2+1)(q^3+1)$. Then
  $M$ has the rank 3 and  $\D_s=\{W\in\O \di \dim(W\cap U)=2\},$ where $|\D_s|=q(q^2+1)(q^2+q+1).$

Noting that  for any two distinct $W_1, W_2\in \O$,  we have $\dim(W_1\cap W_2)=0$ or $2$,   we  introduce the bipartite graph $Y$ with biparts $\O$ and the set $\Psi $ of 2-dim t.s. subspaces of $V$,
where $S\in \Psi$ is adjacent with $W\in \O$ if and only if $S\subset W$. Then  $Y$ has the diameter $d(Y)=4$ and  the girth $g(Y)\ne 4$. By $Y_l(\a)$ we denote the  set of vertices  with the distance $l$ to $\a$.
We  first derive two  statements.

\begin{lem}\label{AB} With the above notations, we have

 (1)\, Let  $W\in \O$ and $S_1$ and $S_2$  be two 2-dim t.s. subspaces of $W$  such that $\dim (S_1\cap S_2)=1$. Then for any point $W_{1i}$ in  $Y_1(S_1)\setminus \{W\}$,  there exists a H-path $W_{1i}\to W_{21}$ of $\G(Y_1(S_1)\cup Y_1(S_2)\setminus \{W\})$.
\vskip 2mm
 (2)\, Let  $S_1$ and $S_2$  be two 2-dim t.s. subspaces such that $S_1+S_2$ be a 4-dim orthogonal subspace. Then there exists a matching between $\G(Y_1(S_1))$ and $\G(Y_1(S_2)).$
\end{lem}
\demo  (1)  Set   $W=\lg u_1',v_1',w_1',x_1'\rg $, $S_1=\lg u_1',w_1' \rg$,  $S_2=\lg u_1',x_1'\rg$  and $V=\lg u_1',u_2'\rg\perp \lg v_1',v_2'\rg\perp \lg w_1',w_2'\rg \perp \lg x_1',x_2'\rg$ be a  decomposition  of   hyperbolic planes for $V$.
 Since
 $$\begin{array}{llll}  Y_1(S_1)\setminus \{W\}&=&\{W_{1k} \di k\in \FF_q\},   \quad  &W_{1k}=\lg u_1',w_1',v_2'+kx_1',kv_1'-x_2' \rg  ,\\
 Y_1(S_2)\setminus \{W\}&=&\{ W_{2k}\di k\in \FF_q\},   &W_{2k}=\lg u_1',x_1',v_2'+kw_1',kv_1'-w_2'\rg ,
\end{array}$$
the point  $W_{20}$ in $Y_1(S_2)\setminus \{W\}$   is adjacent to all points of   $Y_1(S_1)\setminus \{W\}$ and
   $\G(Y_1(S_i))$ is complete for  $i=1,2$,  we may   first pick    a H-path  $W_{1i}\to W_{1j}$  ($i\ne j$) of  $\G(Y_1(S_1)\setminus \{W\})$, and a H-path  $W_{20}\to W_{21}$  of  $\G(Y_1(S_2)\setminus \{W\})$,
  and then get a
 H-path of  $\G(Y_1(S_1)\cup Y_1(S_2)\setminus \{W\})$:  $W_{1i}\to W_{1j}, W_{20}\to W_{21}.$
\vskip 2mm
(2)  Set $S_1=\lg u_1',v_1'\rg$, $S_2=\lg u_2',v_2'\rg$, $W=\lg u_1',v_1',w_1',x_1'\rg\in \O $ and $V=\lg u_1',u_2'\rg\perp \lg v_1',v_2'\rg\perp \lg w_1',w_2'\rg \perp \lg x_1',x_2'\rg$. Then $\G(Y_1(S_1))$ and $\G(Y_1(S_2))$ have a matching
$\lg u_1',v_1',$ $w_1',x_1'\rg\thicksim \lg u_2',v_2',w_1',x_1'\rg$ and
$\lg u_1',v_1',w_2'+kx_1',kw_1'-x_2\rg\thicksim \lg u_2',v_2',w_2'+kx_1',kw_1'-x_2\rg$, where $k\in \FF_q$.\qed

\begin{lem}\label{O2}
Let $M=\O^{+}(8,q)$, acting on  an   orbit $\O$ of maximal totally  singular  subspaces. Then $\G(\O, \Del_s)$  has a H-cycle.
\end{lem}
\demo  Take $U_\infty=U=\lg u_1,v_1,w_1,x_1\rg \in \O$ and $S_{\infty 0}=\lg u_1, v_1\rg \in  \Psi .$ Then
\vspace{-2mm}$$\O_{\infty 0}:=Y_1(S_{\infty 0})=\{ U_\infty, U_k\di k\in\FF_q\}, \, {\rm  where}\,  U_k=\{ \lg u_1, v_1\rg \perp \lg  w_2+kx_1, kw_1-x_2\rg \},\vspace{-2mm}$$
so that $|\O_{\infty 0}|=q+1$.
Noting that $U_\infty\setminus \{S_{\infty 0}\}$  contains  $\frac{(q^4-1)(q^4-q)}{(q^2-1)(q^2-q)}-1=q^4+q^3+2q^2+q$  2-dim t.s.  subspaces, say $S_{\infty j}$, where $1\le j\le q^4+q^3+2q^2+q$, we  set
\vspace{-2mm}$$\O_{\infty j}=Y_1(S_{\infty j})\setminus \{U_\infty\},\,\,{\rm where}\,\,|\O_{\infty j}|=q.\vspace{-2mm}$$
 Noting that for any $k\in \FF_q$,  $U_k$ contains totally $\frac{(q^4-q^2)(q^4-q^3)}{(q^2-1)(q^2-q)}=q^4$  complements of $S_{\infty 0}$, say $S_{ki}$ where $1\le i\le q^4$,  we set
 \vspace{-2mm}$$\O_{ki}=Y_1(S_{ki})\setminus \{U_k\}, \,  {\rm where}\, \,|\O_{ki}|=q.\vspace{-2mm}$$

{\bf Step 1:} {\it Show that   $\O$ is the disjoint union of all subsets  $\O_{\infty j}$ and $\O_{ki}.$}
 \vskip 1.5mm
 Since
\vspace{-2mm}$$|\O_{\infty 0}|+\sum_{j=1}^{q^4+q^3+2q^2+q} |\O_{\infty j}|+\sum_{k\in \FF_q} \sum_{i=1}^{q^4} |\O_{ki}|=q^6+q^5+q^4+2q^3+q^2+q+1=|\O|,\vspace{-2mm}$$
 it suffices to show  all the subset  $\O_{\infty j}$  and $\O_{ki}$ are mutually disjoint, for all  $i, j, k$.

Since  $g(Y)\ne 4$  and $\O_{\infty j}\subset  Y_2(U_{\infty})$, all the subset
   $\O_{\infty j}$ $(0\le j\le q^4+q^3+2q^2+q)$  are mutually  disjoint;  and  for  every given $k$,  all the subset  $\O_{ki}$ ($1\le i\le q^4 $)  are mutually  disjoint.
   Since  $\O_{ki}\subset Y_4(U_{\infty})$,     $\cup _{k,i} \O_{ki}$ and $\cup _j \O_{\infty j}$
are  disjoint.

Suppose that  $W\in \O_{k_1i}\cap \O_{k_2l}.$ Then $W$ contains two 2-dim subspaces:
\vspace{-2mm}$$\begin{array}{lll} S_{k_1i}&=&\lg i_1u_1+i_2v_1+w_2+k_1x_1,i_1'u_1+i_2'v_1+x_2-k_1w_1\rg, \\
S_{k_2l}&=&\lg l_1u_1+l_2v_1+w_2+k_2x_1,l_1'u_1+l_2'v_1+x_2-k_2w_1\rg ,
\end{array}\vspace{-2mm}$$
\f for some  $i_1,i_1',i_2,i_2',l_1,l_1',l_2,l_2'\in \FF_q$.
  Since  $W$ is t.s.,  we get $k_1=k_2$ so that  $S_{k_1i}=S_{k_2l}$ and $\O_{k_1i}=\O_{k_2l}$
(as $g(Y)\ne 4$). Therefore,  all $\O_{\infty j}$ and $\O_{ki}$  are mutually disjoint.

\vskip 1.5mm

{\bf Step 2:} {\it Find a sequence of all 2-dim t.s. subspaces  $S_{\infty j}$ ($j\ne 0$)  (resp. $S_{ki}$) of $U_\infty $  (resp. $U_k$) such that $\dim(S_{\infty j}\cap S_{\infty (j+1)})=1$ (resp. $\dim(S_{ki}\cap S_{k(i+1)})=1$).}
\vskip 1.5mm

Let $\a=\{\a_l \di 0\le l \le q^2-1\}$, $\b=\{\b_l \di 0\le l \le q^2-1\}$ be all vectors of $\lg u_1, v_1\rg $ and $\FF_q=\{ k_l\di 0\le l \le q-1\}$.

   \vskip 1.5 mm
   (a)  $U_\infty$: the set of all  its 2-dim t.s spaces  $S_{\infty j}$ ($j\ne 0$)   is the disjoint union of the following subsets:
   \vspace{-2mm}$$\begin{array}{lll}
    \Lambda(v_1)&=&\{ \lg v_1, y  \rg \di  y \in \lg u_1, x_1, w_1\rg \setminus \lg u_1 \rg \};\\
      \Lambda (u_1+kv_1)&=&\{ \lg u_1+kv_1, y  \rg \di  y \in \lg v_1, x_1, w_1\rg \setminus \lg v_1 \rg \}, \, \, {\rm  where}\, \,  k\in \FF_q;\\
    \Lambda (w_1+\a )&=&\{ \lg w_1+\a, x_1+\b \rg \di  \b \in \lg u_1, v_1 \rg \}, \, \, {\rm  where}\, \, \a\in\lg u_1,v_1\rg.
   \end{array}\vspace{-2mm}$$

  First suppose that $q$ is even.   Take  any sequence:
\vspace{-2mm}$$\begin{array}{ll}
\lg v_1, x_1\rg , \cdots, \lg v_1 ,w_1\rg \,\,{\rm in }\,\, \Lambda(v_1);\\
\lg u_1+k_0 v_1, w_1\rg ,\cdots ,\lg u_1+k_0v_1, x_1\rg \,\, {\rm in}\,\, \Lambda(u_1+k_0v_1);\\
\lg u_1+k_1 v_1, x_1\rg ,\cdots ,\lg u_1+k_1v_1, w_1\rg \,\, {\rm in}\,\, \Lambda(u_1+k_1v_1);\\
\cdots;\\
\lg u_1+k_{q-1}v_1, x_1\rg , \cdots , \lg u_1+k_{q-1}v_1, w_1\rg \,\, {\rm in}\,\, \Lambda(u_1+k_{q-1}v_1);\\
\lg w_1+\a_0, x_1+\a_0\rg , \cdots , \lg w_1+\a_0, x_1+\a_{q^2-1}\rg \,\, {\rm in}\,\, \Lambda(w_1+\a_0),\,\,{\rm where} \,\, \a_0=0;   \\
\lg w_1+\a_1, x_1+\a_{q^2-1}\rg , \cdots , \lg w_1+\a_1, x_1+\a_0 \rg \,\, {\rm in} \,\, \Lambda(w_1+\a_1);\\
\cdots;\\
\lg w_1+\a_{q^2-1}, x_1+\a_{q^2-1}\rg , \cdots ,  \lg w_1+\a_{q^2-1}, x_1+\a_0 \rg \,\, {\rm in}\,\, \Lambda(w_1+\a_{q^2-1}).
\end{array}\vspace{-2mm}$$
Attaching  these sequences, we get a desired  sequence   for all elements   of  $U_\infty\setminus \{S_{\infty 0}\}$, starting with $\lg v_1, x_1\rg $ and ending with   $\lg w_1+\a_{q^2-1}, x_1+\a_0 \rg $.

Second,  suppose that $q$ is odd. With the similar arguments,  we get
a desired  sequence  for all elements   of     $U_\infty\setminus \{S_{\infty 0}\}$:  $\lg v_1, x_1\rg , \cdots , $ $\lg u_1+k_{q-1}v_1, x_1\rg $, $\lg w_1+\a_0, x_1+\a_0\rg $, $\cdots ,$  $\lg w_1+\a_{q^2-1}, x_1+\a_{q^2-1}\rg $.

 \vskip 1.5 mm
   (b)   $U_k$:  recall $U_k=\lg u_1, v_1\rg \perp\lg  w_2+kx_1, x_2-kw_1\rg $
   the set of all  its  complements   $S_{ki}$  of $S_{\infty 0}$  is the disjoint union of the following subset,
    \vspace{-2mm}$$\Lambda ((w_2+kx_1)+\a )=\{ \lg (w_2+kx_1)+\a,  (x_2-kw_1)+\b \rg \di  \b \in \lg u_1, v_1 \rg \}, \, \,  {\rm where}\,\,  \a\in \lg u_1, v_1\rg .\vspace{-2mm}$$
   By replacing $w_1$ and $x_1$   in the third part $\Lambda (w_1+\a )$ of the  case (a)  by $w_2+kx_1$ and $x_2-kw_1$ respectively,    we get a desired full  sequence:
    $\lg w_2+kx_1+\a_0,x_2-kw_1+\a_0\rg \to $ either $\lg w_2+kx_1+\a_{q^2-1},x_2-kw_1+\a_{0}\rg$ if $q$ is even or  $\lg w_2+kx_1+\a_{q^2-1},x_2-kw_1+\a_{q^2-1}\rg$, if $q$ is odd.

\vskip 1.5mm
{\bf Step 3:} {\it Find a  H-cycle for $\G$.}
\vskip 1.5mm
 Set
$$\G_\infty:=\G(\cup _{j=1}^{q^4+q^3+2q^2+q} \O_{\infty j})\, \, {\rm and}\, \,\G_{k_i}=\G(\cup _{j=1}^{q^4} \O_{k_i j}),$$
 where $\FF_q=\{k_i \di 0\le i \le q-1\}$. By
$W_{\infty j; k}(1\le j\le q^4+q^3+2q^2+q$)  and $W_{k_i j; k}$  ($1\le j\le q^4$), we denote  points in $\O_{\infty j}$ and  $\O_{k_i j}$, respectively.

 Since $S_{\infty,q^4+q^3+2q^2+q}+S_{k_0,1}$  (resp. $S_{k_i, q^4}+S_{k_{i+1},1}$   ($0\le i\le q-1$)) is a 4-dim orthogonal space,
  by Statement B, there exists a matching between $\O_{\infty,q^4+q^3+2q^2+q}$ and $\O_{k_0,1}$  (resp. $\O_{k_i, q^4}$ and $\O_{k_{i+1},1}$).
  Assume that the point  $W_{\infty ,  q^4+q^3+2q^2+q;  1}$  in  $\O_{\infty,q^4+q^3+2q^2+q}$  is adjacent to the point $W_{k_0,1; l_0}$ in  $\O_{k_0,1}$   for some $l_0$; and
  the point  $W_{k_i, q^4; 1}$  in  $\O_{k_i, q^4}$  is adjacent to the point $W_{k_{i+1},1; l_{i+1}}$ in  $\O_{k_{i+1}, 1}$.

First consider the induced subgraph  $\G_\infty$.
  Recall $\O_{\infty j}=Y_1(S_{\infty j})\setminus \{U_\infty\}$, for any $j$.   Then
   $\O_{\infty j}$ consists of  all  $W_{\infty j; k}=\lg S_{\infty j},  v_2'+kx_1',kv_1'-x_2' \rg ,$  for some $x_1', x_2', v_1', v_2'\in V$,  where
     $k\in \FF_q$.
In {\bf Step 2}, we have constructed   a sequence of all $S_{\infty j}$  where $1\le j\le q^4+q^3+2q^2+q$,
 starting with $S_{\infty 0}=\lg v_1, x_1\rg $ and ending with   either $S_{\infty (q^4+q^3+2q^2+q)}=\lg w_1+\a_{q^2-1}, x_1+\a_0 \rg $ if $q$ is even or   $\lg w_1+\a_{q^2-1}, x_1+\a_{q^2-1}\rg $ if $q$ is odd.
    Along this sequence and by Statement A,   we get  a H-path  $H_\infty :$ $W_{\infty 1;0}\to W_{\infty  (q^4+q^3+2q^2+q);1}$
     of $\G_\infty $.

 Second,  consider the induced subgraph  $\G_{k_i}$, where $0\le i\le q-1$.    Recall $\O_{k_i j}=Y_1(S_{k_i j})\setminus \{U_{k_i}\}$.   Then
   $\O_{k_i j}$ consists of  all  $W_{k_i j; k}=\lg S_{k_i j}, w_2'+kx_1',kw_1'-x_2\rg $ for some $x_1', x_2', w_1', w_2'\in V$ and
     $k\in \FF_q$.
In {\bf Step 2}, we have constructed   a sequence of all $S_{k_i j}$  ($1\le j\le q^4$),
 starting with $\lg w_2+kx_1+\a_0,x_2-kw_1+\a_0\rg$ and ending with  either $\lg w_2+kx_1+\a_{q^2-1},x_2-kw_1+\a_{0}\rg$ if $q$ is even or  $\lg w_2+kx_1+\a_{q^2-1},x_2-kw_1+\a_{q^2-1}\rg$, if $q$ is odd.
 Along this sequence and by Statement A,  we get  the corresponding H-path  $H_{k_i}:$  $W_{k_i 1; l_i}\to W_{k_i q^4; 1}.$
     of $\G_{k_i}$.

From the second paragraph, the ending point of  $H_\infty$ (resp. $H_{k_i}$) is adjacent with the starting point of  $H_{k_0}$ (resp. $H_{k_{i+1}}$).
Thus we get   a H-cycle $C$ of $\G$ (see Figure 7.)
$$C:\quad U_\infty, H_\infty, H_{k_0},\cdots,H_{k_{q-1}},U_{k_{q-1}},\cdots, U_\infty.\vspace{-2mm}$$
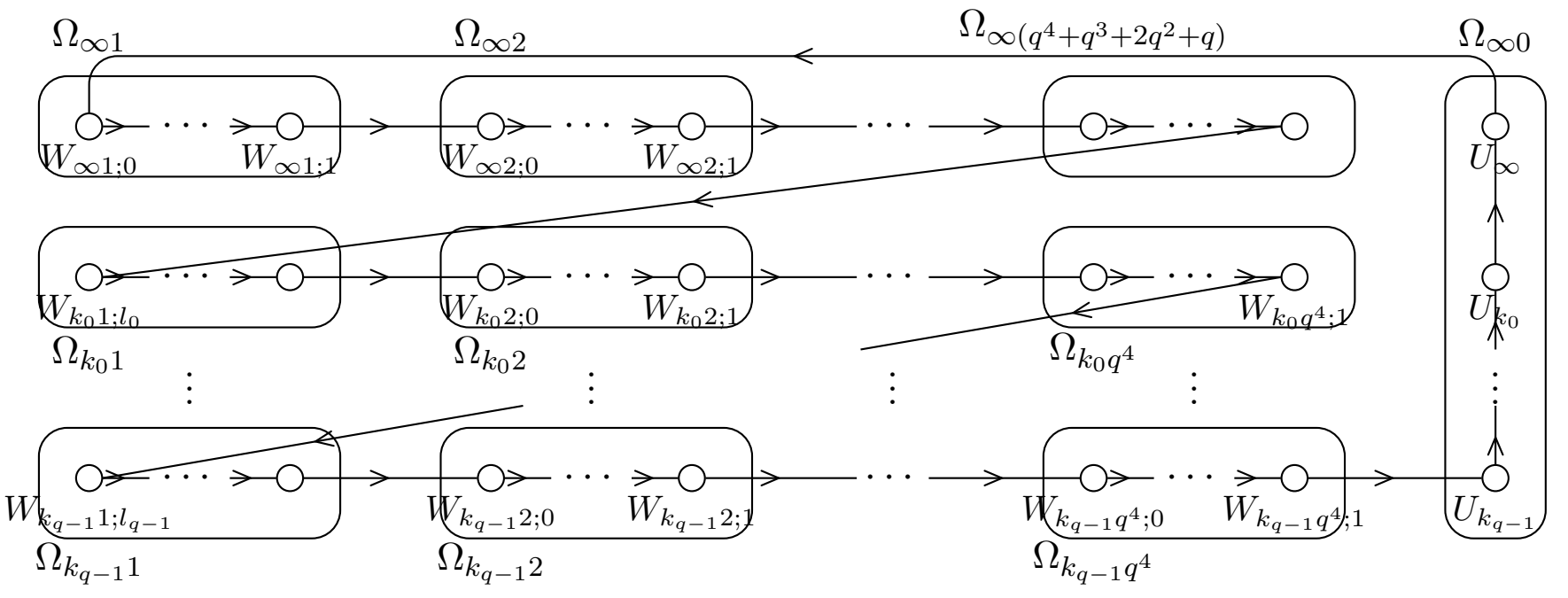
\begin{figure}[htbp]
  \centering
  \begin{tikzpicture}[
    scale=0.75, transform shape,
    p/.style={draw, circle, fill=white, inner sep=0pt, minimum size=2.6mm},
    d/.style={fill=white, inner sep=2pt},
    box/.style={draw, rounded corners=6pt},
    lab/.style={anchor=north, yshift=3pt, font=\footnotesize},
    tlab/.style={anchor=south, font=\small},
    midarrow/.style={
      postaction=decorate,
      decoration={markings, mark=at position 0.5 with {\arrow{angle 45}}}
    },
    remidarrow/.style={
      postaction=decorate,
      decoration={markings, mark=at position 0.5 with {\arrowreversed{angle 45}}}
    }
  ]
    \node[p] (a10)  at (0, 0)  {};
    \node[d] (d10)  at (1, 0)  {$\cdots$};
    \node[p] (a20)  at (2, 0)  {};
    \draw[midarrow] (a10.east) -- ($(a10)!0.3!(a20)$);
    \draw[midarrow] ($(a10)!0.7!(a20)$) -- (a20.west);
    \node[p] (a30)  at (4, 0)  {};
    \node[d] (d20)  at (5,0)  {$\cdots$};
    \node[p] (a40)  at (6, 0)  {};
    \draw[midarrow] (a30.east) -- ($(a30)!0.3!(a40)$);
    \draw[midarrow] ($(a30)!0.7!(a40)$) -- (a40.west);
    \node[d] (d30)  at (8,0)  {$\cdots$};
    \node[p] (a50)  at (10, 0)  {};
    \node[d] (d40)  at (11, 0)  {$\cdots$};
    \node[p] (a60)  at (12, 0)  {};
    \draw[midarrow] (a50.east) -- ($(a50)!0.3!(a60)$);
    \draw[midarrow] ($(a50)!0.7!(a60)$) -- (a60.west);
    \node[p] (a70)  at (14, 0)  {};
    \draw[midarrow] (a20.east) -- (a30.west);
    \draw[midarrow] (a40.east) -- (d30.west);
    \draw[midarrow] (d30.east) -- (a50.west);

    \node[p] (a11)  at (0, -1.5)  {};
    \node[d] (d11)  at (1, -1.5)  {$\cdots$};
    \node[p] (a21)  at (2, -1.5)  {};
    \draw[midarrow] (a11.east) -- ($(a11)!0.3!(a21)$);
    \draw[midarrow] ($(a11)!0.7!(a21)$) -- (a21.west);
    \node[p] (a31)  at (4, -1.5)  {};
    \node[d] (d21)  at (5, -1.5)  {$\cdots$};
    \node[p] (a41)  at (6, -1.5)  {};
    \draw[midarrow] (a31.east) -- ($(a31)!0.3!(a41)$);
    \draw[midarrow] ($(a31)!0.7!(a41)$) -- (a41.west);
    \node[d] (d31)  at (8,-1.5)  {$\cdots$};

    \node[p] (a51)  at (10, -1.5)  {};
    \node[d] (d41)  at (11, -1.5)  {$\cdots$};
    \node[p] (a61)  at (12, -1.5)  {};
    \draw[midarrow] (a51.east) -- ($(a51)!0.3!(a61)$);
    \draw[midarrow] ($(a51)!0.7!(a61)$) -- (a61.west);
    \node[p] (a71)  at (14, -1.5)  {};
    \draw[midarrow] (a21.east) -- (a31.west);
    \draw[midarrow] (a41.east) -- (d31.west);
    \draw[midarrow] (d31.east) -- (a51.west);

    \node[d] (v1)  at (1, -2.5)  {$\vdots$};
    \node[d] (v1)  at (5, -2.5)  {$\vdots$};
    \node[d] (v1)  at (8, -2.5)  {$\vdots$};
    \node[d] (v1)  at (11, -2.5)  {$\vdots$};
    \node[d] (v1)  at (14, -2.5)  {$\vdots$};
    \node[p] (a12)  at (0, -3.5)  {};
    \node[d] (d12)  at (1, -3.5)  {$\cdots$};
    \node[p] (a22)  at (2, -3.5)  {};
    \draw[midarrow] (a12.east) -- ($(a12)!0.3!(a22)$);
    \draw[midarrow] ($(a12)!0.7!(a22)$) -- (a22.west);
    \node[p] (a32)  at (4, -3.5)  {};
    \node[d] (d22)  at (5, -3.5)  {$\cdots$};
    \node[p] (a42)  at (6, -3.5)  {};
    \draw[midarrow] (a32.east) -- ($(a32)!0.3!(a42)$);
    \draw[midarrow] ($(a32)!0.7!(a42)$) -- (a42.west);

    \node[d] (d32)  at (8,-3.5)  {$\cdots$};

    \node[p] (a52)  at (10, -3.5)  {};
    \node[d] (d42)  at (11, -3.5)  {$\cdots$};
    \node[p] (a62)  at (12, -3.5)  {};
    \draw[midarrow] (a52.east) -- ($(a52)!0.3!(a62)$);
    \draw[midarrow] ($(a52)!0.7!(a62)$) -- (a62.west);
    \node[p] (a72)  at (14, -3.5)  {};
    \draw[midarrow] (a22.east) -- (a32.west);
    \draw[midarrow] (a42.east) -- (d32.west);
    \draw[midarrow] (d32.east) -- (a52.west);

    \draw[midarrow] (a60.west) -- (a11.east);
    \draw[midarrow] (a61.west) -- ($(a61)!0.36!(a12)$);
    \draw[midarrow] ($(a61)!0.64!(a12)$) -- (a12.east);
    \draw[midarrow] (a62.east) -- (a72.west);
    \draw[midarrow] (a72.north) -- ($(a72)!0.36!(a71)$);
    \draw[midarrow] ($(a72)!0.64!(a71)$) -- (a71.south);
    \draw[midarrow] (a71.north) -- (a70.south);
    \draw[midarrow, rounded corners=6pt]
    (a70.north) -- (14, 0.7) -- (0,0.7) -- (a10.north);

    \draw[box] (-0.5,0.5) rectangle (2.5, -0.5);
    \draw[box] (-0.5,-1) rectangle (2.5, -2);
    \draw[box] (-0.5,-3) rectangle (2.5, -4.1);
    \draw[box] (3.5,0.5) rectangle (6.6, -0.5);
    \draw[box] (3.5,-1) rectangle (6.6, -2);
    \draw[box] (3.5,-3) rectangle (6.6, -4.1);
    \draw[box] (9.5,0.5) rectangle (12.6, -0.5);
    \draw[box] (9.5,-1) rectangle (12.6, -2);
    \draw[box] (9.5,-3) rectangle (12.5, -4.1);
    \draw[box] (13.5,0.5) rectangle (14.5, -4.1);

    \node[lab] at (a70.south) {$U_\infty$};
    \node[lab] at (a71.south) {$U_{k_0}$};
    \node[lab] at (a72.south) {$U_{k_{q-1}}$};
    \node[lab] at (a10.south) {$W_{\infty1;0}$};
    \node[lab] at (a20.south) {$W_{\infty1;1}$};
    \node[lab] at (a30.south) {$W_{\infty2;0}$};
    \node[lab] at (a40.south) {$W_{\infty2;1}$};
    \node[lab] at (a11.south) {$W_{k_0 1;l_0}$};
    \node[lab] at (a31.south) {$W_{k_0 2;0}$};
    \node[lab] at (a41.south) {$W_{k_0 2;1}$};
    \node[lab] at (a61.south) {$W_{k_0 q^4;1}$};
    \node[lab] at (a12.south) {$W_{k_{q-1} 1;l_{q-1}}$};
    \node[lab] at (a32.south) {$W_{k_{q-1} 2;0}$};
    \node[lab] at (a42.south) {$W_{k_{q-1} 2;1}$};
    \node[lab] at (a52.south) {$W_{k_{q-1} q^4;0}$};
    \node[lab] at (a62.south) {$W_{k_{q-1} q^4;1}$};

    \node[tlab] at (0,0.6) {$\O_{\infty 1}$};
    \node[tlab] at (0,-2.6) {$\O_{k_0 1}$};
    \node[tlab] at (0,-4.7) {$\O_{k_{q-1} 1}$};

    \node[tlab] at (4,0.6) {$\O_{\infty 2}$};
    \node[tlab] at (4,-2.6) {$\O_{k_0 2}$};
    \node[tlab] at (4,-4.7) {$\O_{k_{q-1} 2}$};

    \node[tlab] at (10,0.6) {$\O_{\infty(q^4+q^3+2q^2+q)}$};
    \node[tlab] at (10,-2.6) {$\O_{k_0 q^4}$};
    \node[tlab] at (10,-4.7) {$\O_{k_{q-1} q^4}$};

    \node[tlab] at (14,0.6) {$\O_{\infty 0}$};
 \end{tikzpicture}
  \caption{{\footnotesize H-cycle of the graph $\G$ in Lemma 6.1.}}
  \vspace{-2mm}
\end{figure}

\subsection{$M=\O^{+}(10,q)$}
\begin{lem}\label{O3}
Let $M=\O^{+}(10,q)$, acting on  an $M$-orbit $\O$ of maximal totally  singular  subspaces. Then $\G(\O, \Del_s)$  has a H-cycle.
\end{lem}
\demo
Let $V=\lg u_1,u_2\rg\perp\lg v_1,v_2\rg\perp\lg w_1,w_2\rg\perp\lg x_1,x_2\rg\perp\lg y_1,y_2\rg$ be a decomposition of hyperbolic planes for $V$. Let $U=\lg u_1,v_1,w_1,x_1,y_2\rg$ be an m.t.s. subspace of $V$. The group $\GO^+(10,q)$ is transitive on the set of m.t.s. subspaces of $V$ but $M$ has two orbits, say $\O$ and $\O'$, while $\O$ (resp. $\O'$) consists of all those intersecting $U$ with odd (resp. even) dimension. Since the actions of $M$ on $\O$ and $\O'$ are equivalent, we just consider its action on $\O$, where $|\O|=(q+1)(q^2+1)(q^3+1)(q^4+1)$. Then $\D_s=\{W\in\O \di \dim(W\cap U)=3\}$, where $|\D_s|=q(q^2+1)(q^4+q^3+q^2+q+1)$.

Noting that for any two distinct $W_1,W_2\in \O$, we have $\dim(W_1\cap W_2)=1$ or 3, we introduce the bipartite $Y$ with biparts $\O$ and the set $\Psi$ of 3-dim t.s. subspaces of $V$, where $S\in\Psi$ is adjacent with $W\in \O$ if and only if $S\subset W$.
 Similar to the case (1), $d(Y)=4$ and $g(Y)\ne 4$. Also one has
\vskip 1.5mm
{\it Statement A: Let  $W\in \O$ and $S_1$ and $S_2$  two 3-dim t.s. subspaces of $W$  such that $\dim (S_1\cap S_2)=2$.
Then for any point $W_{1i}$ in $Y_1(S_1)\setminus \{W\}$, there exists a H-path $W_{1i}\to W_{21}$ of $\G(Y_1(S_1)\cup Y_1(S_2)\setminus \{W\})$.}
\vskip 1.5mm
Let $U_\infty=U\in\O$. Then $U_\infty$ contains $m=\frac{(q^5-1)(q^5-q)(q^5-q^2)}{(q^3-1)(q^3-q)(q^3-q^2)}=(q^2+1)(q^4+q^3+q^2+q+1)$ 3-dim t.s. subspaces, say $S_{\infty j}$, where $1\le j \le m$. Set
\vspace{-2mm}$$\O_{\infty j}=Y_1(S_{\infty j})\setminus \{U_\infty\},\,\, {\rm where}\,\, |\O_{\infty j}|=q.\vspace{-2mm}$$

Now, we introduce another  bipartite  graph $Z$ with  biparts the set $\PG(U)$ and the set $\Psi$ of 3-dim t.s. subspaces of $U$, where $\lg v\rg\in \PG(U)$ is adjacent with $S\in \Psi$ if and only if $v\in S$.
For any subset $T\subset \PG(U)$, let $Z_1(T)=\cup _{\lg v\rg\in T} Z_1(\lg v\rg ).$
Since  $|Z_1(\lg v\rg)|=(q^2+1)(q^2+q+1)$ and $|Z_1(S)|=q^2+q+1$, we have $|T|\le|Z_1(T)|$. By P. Hall's Theorem \cite[Theorem 1.1.3]{LP86}, the bipartite graph $Z$ has a matching $\phi:$ $\lg v\rg \in \PG(U)\to \phi(\lg v\rg)\in \Psi$.
 Then for $\lg v\rg\in \PG(U)$,   we have $\phi(\lg v\rg )=S_{\infty j}$ for some $j$ so that  the set  $\O_{\infty j}= Y_1(\phi(\lg v\rg ))\setminus U_\infty$ can be expressed by $\{ U_k(\lg v\rg )\di k\in \FF_q\}$.

Note that $U_k(\lg v\rg )$ contains totally $\frac{(q^4-q^2)(q^4-q^3)}{(q^2-1)(q^2-q)}=q^4$ 3-dim t.s. subspaces intersecting $\phi(\lg v\rg )$ in $\lg v\rg$, say $S_{ki}(\lg v\rg )$, where $1 \le i\le q^4$, we set
\vspace{-2mm}$$\O_{ki}(\lg v\rg )=Y_1(S_{ki}(\lg v\rg ))\setminus U_k(\lg v\rg ),\,\, {\rm where}\,\, |\O_{ki}(\lg v\rg )|=q.\vspace{-2mm}$$
According to Lemma~\ref{O2}, we can observe that $U_k(\lg v\rg)$ is equivalent to $U_k$, $S_{ki}(\lg v\rg)$ is equivalent to $S_{ki}$ so that $\O_{ki}(\lg v\rg)$ is equivalent to $\O_{ki}$, where $k\in \FF_q$, $1\le i \le q^4$.

{\it Example: Take $\lg y_1\rg \subset U$ and $\phi(\lg y_1\rg)=\lg u_1,v_1,y_1\rg$. Then $U_k(\lg y_1\rg)=\{\lg u_1,v_1,y_1,$ $w_2+kx_1,x_2-kw_1\rg \di k\in\FF_q\}$, and $U_k=\{\lg u_1,v_1,w_2+kx_1,x_2-kw_1\rg \di k\in\FF_q\}$ in Lemma~\ref{O2}. We can observe that $U_k(\lg y_1\rg)/\lg y_1\rg$ corresponds to $U_k$, $S_{ki}(\lg y_i\rg)$ intersecting $\lg u_1,v_1,y_1\rg$ in $\lg y_1\rg$ corresponds to $S_{ki}$ intersecting $\lg u_1,v_1\rg$ trivially, and $\O_{ki}(\lg y_1\rg)$ corresponds to $\O_{ki}$ by removing $\lg y_1\rg$.}

Then we may get a  H-path $H_k(\lg v\rg )$ in $\G(\cup_{i=1}^{q^4}\O_{ki}(\lg v\rg ))$ corresponding to $H_k$ in $\G(\cup_{i=1}^{q^4}\O_{ki})$ by same method in Lemma~\ref{O2}.
Moreover, we can have a H-path in $\G(\O_{\infty j}\cup (\cup_{k,i}\O_{ki}(\lg v\rg )))$,
\vspace{-2mm} $$ L(\lg v\rg )=U_{k_0}(\lg v\rg ),H_{k_0}(\lg v\rg ),H_{k_1}(\lg v\rg ),\cdots,H_{k_{q-1}}(\lg v\rg ),U_{k_{q-1}}(\lg v\rg ),\cdots,U_{k_1}(\lg v\rg ),\vspace{-2mm}$$
corresponding to a path in Lemma~\ref{O2}: $U_{k_0},H_{k_0},H_{k_1},\cdots,H_{k_{q-1}},U_{k_{q-1}},\cdots,U_{k_2},U_{k_1}$ by removing $U_\infty$ and $H_\infty$ for H-cycle.

{\bf Step 1:} {\it Show that $\O$ is the disjoint union of all subset $\O_{\infty j},\O_{ki}(\lg v\rg )$, where $\lg v\rg \in \PG(U)$.}
\vskip 1.5mm
Since
$$1+\sum_{j=1}^{m}|\O_{\infty j}|+\sum_{\lg v\rg\in \PG(U)}\sum_{k\in\FF_q} \sum_{i=1}^{q^4}|\O_{ki}(\lg v\rg)|=|\O|,$$
it suffices to show all $\O_{\infty j}$ and $\O_{ki}(\lg v\rg )$ are mutually disjoint, for any $i,j,k,v$.
Since $g(Y)\ne 4$, $\O_{\infty j}\subset Y_2(U_\infty)$ and $\O_{ki}(\lg v\rg )\subset Y_4(U_\infty)$, according to the {\bf Step 1} in Lemma~\ref{O2}, $\O_{\infty j}$ and $\O_{ki}(\lg v\rg )$ are mutually disjoint for any $i,j,k$. For $\lg u\rg\ne \lg v\rg$, clearly $\O_{ki}(\lg u\rg )$ and $\O_{ki}(\lg v\rg )$ are all m.s.t. spaces intersecting $U_\infty$ in $\lg u\rg$ and $\lg v\rg$ respectively, then they are mutually disjoint.
\vskip 1.5mm

{\bf Step 2:} {\it Find a sequence of all 3-dim t.s. spaces $S_{\infty j}$ of $U_\infty$ such that $\dim(S_{\infty j}\cap S_{\infty (j+1)})=2$.}
\vskip 1.5mm

In 3-dim subspace $\lg u_1,v_1,w_1\rg\subset U_\infty$, let $\{\lg h_l\rg \di 1\le l \le q^2+q+1\}$ be all 1-dim subspaces so that $H_l$ be the complement of $\lg h_l\rg$. Similarly, we have a bipartite graph with biparts $\{\lg h_l\rg\}$ and $\{H_l\}$, where $\lg h_l\rg$ is adjacent with a 2-dim subspace, say $H_l'$, if and only if $h_l\in H_l'$. Then we have a matching $\lg h_l\rg\thicksim H_l'$ by P. Hall Theorem, where $1\le l \le q^2+q+1$.
The set of all its 3-dim t.s. subspaces $S_{\infty j}$ of  $U_\infty$ is the disjoint union of the following subset:
\vspace{-2mm}$$\begin{array}{lll}
\{ \lg u_1,v_1,w_1\rg \};\\
\Lambda(H_l')&=&\{ \lg H_l', z \rg \di  z \in \lg h_l', x_1,y_1\rg \setminus \lg h_l' \rg \};\\
\Lambda(h_l,x_1+\a)&=&\{\lg h_l,x_1+\a, y_1+\b \rg \di \b \in H_l\};\\
&&{\rm where}\,\, \a \in H_l,\, 1\le l \le q^2+q+1.
\end{array}\vspace{-2mm}$$

First suppose that $q$ is even. Starting with $\lg u_1,v_1,w_1\rg$, take any sequence:
\vspace{-2mm}$$\begin{array}{ll}
\lg H_1',x_1\rg,\cdots,\lg H_1',x_1+y_1\rg \,\, {\rm in }\,\, \Lambda(H_1')\setminus \lg H_1',y_1\rg;\\
\lg h_1,x_1+\a_0,y_1+\a_0\rg,\cdots,\lg h_1,x_1+\a_{0},y_1+\a_{q^2-1}\rg \,\, {\rm in } \,\, \Lambda(h_1,x_1+\a_0);\\
\lg h_1,x_1+\a_1,y_1+\a_{q^2-1}\rg,\cdots,\lg h_1,x_1+\a_1,y_1+\a_0\rg \,\, {\rm in } \,\, \Lambda(h_1,x_1+\a_1);\\
\cdots;\\
\lg h_1,x_1+\a_{q^2-1},y_1+\a_{q^2-1}\rg,\cdots,\lg h_1,x_1+\a_{q^2-1},y_1+\a_0\rg \,\, {\rm in } \,\, \Lambda(h_1,x_1+\a_{q^2-1});\\
\lg H_1',y_1\rg;\quad \lg H_2',y_1\rg,\cdots,\lg H_2',x_1+y_1\rg \,\, {\rm in }\,\, \Lambda(H_2')\setminus \lg H_2',x_1\rg;\\
\lg h_2,x_1+\a_0,y_1+\a_0\rg,\cdots,\lg h_2,x_1+\a_{q^2-1},y_1+\a_0\rg;\\
\cdots;\\
\lg h_2,x_1+\a_{q^2-1},y_1+\a_{q^2-1}\rg,\cdots,\lg h_2,x_1+\a_0,y_1+\a_{q^2-1}\rg; \lg H_2',x_1\rg;\\
\cdots;\\
\lg H_{q^2+q+1}',x_1\rg,\cdots,\lg H_{q^2+q+1}',x_1+y_1\rg \,\, {\rm in }\,\, \Lambda(H_{q^2+q+1}')\setminus \lg H_{q^2+q+1}',y_1\rg.\\
\cdots; \lg H_{q^2+q+1}',y_1\rg,\,\,{\rm where\,\, \a_0=0}.
\end{array}\vspace{-2mm}$$
Then we get a desired sequence for all elements of $U_\infty$:  $\lg u_1,v_1,w_1\rg ,  \cdots , \lg H_{q^2+q+1}',y_1\rg$.

Second, suppose that $q$ is odd. With the similar arguments, we get a desired sequence for all elements of $U_\infty$: $\lg u_1,v_1,w_1\rg, \cdots, \lg H_{q^2+q+1}', y_1\rg$.

\vskip 1.5mm
{\bf Step 3:} {\it Find a H-cycle for $\G$.}
\vskip 1.5mm
Set $\G_\infty:=\G(\cup_{j=1}^{m}\O_{\infty j})$, while   points  in $\O_{\infty j}$ are denoted by $W_{\infty j;k}$.
First consider the induced subgraph $\G_\infty$. In {\bf Step 2}, we have constructed a sequence of all $S_{\infty j}$, where $1\le j\le m$: $$S_{\infty 1}=\lg u_1,v_1,w_1\rg, \cdots,S_{\infty m}=\lg H_{q^2+q+1}',y_1\rg.$$  Along this sequence and by Statement A,  we get  the corresponding H-path   $H_\infty$: $W_{\infty 1;0}\to W_{\infty m;1}$ of $\G_\infty$.

Recall that we may find $q^4+q^3+q^2+q+1$ 3-dim t.s. subspaces  in $\{S_{\infty j}\}$, which  corresponds to all points in $\PG(U)$, saying  $\{S_{\infty,j_s}=\phi (\lg v_s\rg )\di 1 \le s \le q^4+q^3+q^2+q+1\}$, where $\lg v_s\rg\in\PG(U)$.
 Replacing  the H-path $W_{\infty,j_s;0}\to W_{\infty,j_s;1}$ of the complete graph $\G(\O_{\infty, j_s})$  by $L(\lg v_s\rg )$ in $H_\infty$, where $1 \le s \le q^4+q^3+q^2+q+1$, and we denote this path by $H$. Finally we get a H-cycle for $\G$:
$U_\infty, H ,U_\infty.$
\qed

\section{Actions on   nonsingular hyperplanes }
The main result of this section is the following lemma.
\begin{lem}\label{7}
Let  $M=\O(2m+1,q)$ where $q\in \{ 4, 8\}$, acting on an orbit $\O$. Then   $\G=\G(\O, \D_s)$ has a H-cycle.
\end{lem}

 By Eq(\ref{odd}), let  $V=\perp_{i=1}^m\lg u_i,v_i\rg  \perp \lg  w\rg $ be an orthogonal space of dimensional $2m+1$,  where
 $$Q(u_i)=Q(v_i)=0, \quad f(u_i,v_i)=1,\quad Q(w)=1, \forall 1\le i\le m.$$
For convince, a   hyperplane $\PG(U)$ of  $\PG(V)$ will be denoted $U$ too.
 A hyperplane $U$ is said to be {\it nonsingular} if   $\rad(Q\di _U)=0$ (so $\rad (f\di_U)=0$).
  Let $\Lambda $ be  the set of  all nonsingular hyperplanes of  $\PG(V)$. Then we have

\begin{lem}\label{+-}  The group  $M=\O(2m+1,q)$ has two orbits on $\Lambda$,  say $ \Lambda^+$ and $ \Lambda^-$, where
 $|\Lambda^+|=\frac{q^{2m}+q^m}{2}$   and $|\Lambda^-|=\frac{q^{2m}-q^m}{2}$ so that  $|\Lambda |=q^{2m}$.
\end{lem}
\demo  Pick up an $U_0=\perp_{i=1}^m\lg u_i,v_i\rg $ in $\Lambda $. Since for any $U\in  \Lambda$, we have ${\overline U}={\overline U_0}$
in the quotient space $V/\lg w\rg $. Therefore $U=\perp_{i=1}^m\lg u_i+\ld_iw,v_i+\mu_iw\rg$, where $\ld_i,\mu_i\in \FF_q$,  so that  $|\Lambda |=q^{2m}$.
   Moreover,  $M=\O(2m+1,q)$ has two orbits on $\Lambda$, saying $ \Lambda^+$ and $ \Lambda^-$, consisting of    hyperplanes   of plus  and minus  type, respectively.

   For  every plane $\lg u, v\rg $, it is  of plus  type if and only if
  it contains  some singular vectors.
Consider all planes  $P=\lg u_i+\ld_iw, v_i+\mu_iw\rg $.
Check that  for any $(x,y)\ne (0,0)$, $P\in \Lambda^+$ if and only if
\vspace{-5pt}   $$Q(x(u_i+\ld_iw)+y(v_i+\mu_iw))=x^2\ld_i^2+y^2\mu_i^2+xy=0,$$
that is $\ld_i\mu_i\in x_1^2+x_1$ for some $x_1\in \FF_q$.
Therefore, given $\lg u_i, v_i\rg $  for $q=8$ (resp $q=4$) we have   $36$  (resp 10)  planes $P$ of plus type, including $U_0$ and then  28 (resp 6) ones of  minus type.
Moreover, $U=\perp_{i=1}^m\lg u_i+\ld_iw,v_i+\mu_iw\rg$ is of plus type if and only if it contains  even number of  minus type planes  $\lg u_i+\ld_iw,v_i+\mu_iw\rg $.
For $q=8$ (resp 4), set $r_q=28$ (resp. 6) and $s_q=36$ (resp. 10). Then we have
\vspace{-5pt}  $$|\Lambda^+|= \sum_{r=0}^{\lfloor \frac{m}{2}\rfloor}\binom{m}{2r}r_q^{2r}s_q^{m-2r}=\frac{q^{2m}+q^m}{2},$$
\f so
$  |\Lambda^-|= |\Lambda|- |\Lambda^+|=\frac{q^{2m}-q^m}{2}.$
\qed
\vskip 2mm
In what follows we discuss $\Lambda^+$ and $\Lambda^-$, separately, so that Lemma~\ref{7} will be proved in by  Lemmas~\ref{+1} and ~\ref{-1}.
\subsection{ $\O=\Lambda ^+$}
Recall  $U_0=\perp_{i=1}^m\lg u_i,v_i\rg$. The we have
\begin{lem}\label{+1} The stabilizer $M_{U_0}$ has two nontrivial orbits on $\O$:
$$\begin{array}{lll}
&\D_1=\{W\in\O \di \dim(\rad (Q\di_{U_0\cap W})=1)\},\quad &|\D_1|=(q^m-1)(q^{m-1}+1);\\

&\D_2=\{W\in\O \di \dim(\rad (Q\di_{U_0\cap W})=0)\}, \quad  &|\D_2|=(\frac q2-1)\cdot q^{m-1}(q^m-1).\end{array}$$
\f If $q=4$, then $|\D_1|, |\D_2|\ge \frac 13|\O|$.\end{lem}
\demo  For  any $W\in \Lambda^+\setminus \{U_0\}$,     we have $\dim(U_0\cap W)=2m-1$, which is odd, so that  $\dim(\rad(f\di_{U_0\cap W}))=1$. Therefore,    $M_{U_0}$ has two nontrivial orbits as shown in the lemma.
For any maximal subspace $H$ of  $U_0$, set $\phi(H)=H^\perp\cap U_0$.   Then $\phi(H)\subset H$.

First suppose  $\rad(Q\di _H)=1$. Then  $Q(\phi(H))=0$.
For instance, take $H=\lg u_1 \rg\perp (\perp_{i=2}^m\lg u_i,v_i\rg )$ and $\phi(H)=\lg u_1\rg $.
Then  in $\Lambda^+\setminus \{U_0\}$, the hyperplanes $W$ of  plus type containing $H$ are of the form  $\lg u_1,v_1+kw\rg \perp (\perp_{i=2}^m\lg u_i,v_i\rg)$, where $k\in \FF_8^*$.
Since there are exactly  $\frac{(q^m-1)(q^{m-1}+1)}{q-1}$ singular points $\phi(H)$,  we  get $|\D_1|=(q^m-1)(q^{m-1}+1)$.

Second,  suppose  $\dim(\rad(Q\di _H))=0$. Then  $Q(\phi(H))=1$. For instance, take
$H=\lg u_1+v_1 \rg\perp (\perp_{i=2}^m\lg u_i,v_i\rg )$ and $\phi(H)=\lg u_1+v_1\rg $.
Then in $\Lambda^+\setminus \{U_0\}$, the  hyperplanes of plus type containing $H$ are $\lg u_1+kw, v_1+kw\rg\perp(\perp_{i=2}^m\lg u_1,v_1\rg)$, where $k=x^2+x$ for some $x\in \FF_8^*$.
Since there are exactly  $q^{m-1}(q^m-1)$ nonsingular points $\phi(H)$ in $\PG(U_0)$,   we have  $|\D_2|=(\frac q2-1)\cdot q^{m-1}(q^m-1)$.

Clearly, if $q=4$, then $|\D_1|, |\D_2|\ge \frac 13|\O|$.\qed

\vskip 2mm
 So from now on, let $q=8$ so that  $\D_s=\D_1$, that is two distinct $W_1$ and $W_2$ in $\O$ are adjacent if and only if $\rad(Q\di _{W_1\cap W_2})=0$.
 For any  maximal subspace $H$ of $U_0$, let $\o H$ be  the set of all  hyperplanes  containing $H$ in $\Lambda^{+}\setminus \{U_0\}$. Then the induced subgraph  $\G(\o{H})$ is either a complete graph $K_7$ or  an empty graph $\o K_3$,
 provided  $\phi(H)$ is either  singular or nonsingular (as $q=8$).
First we  derive a lemma.

\begin{lem}\label{+2}
Let $H_1$ and $H_2$ be  two maximal subspaces of $U_0$ such that $\phi(H_1)$ and $\phi(H_2)$ are two  singular points (resp. nonsingular points) which are mutually orthogonal. Then $\G(\o{H_1} \cup \o{H_2})$ is $K_{14}$ (resp.  $3K_2$).
\end{lem}
\demo
First assume that $\phi(H_1)$ and $\phi(H_2)$ are singular. No loss,  let $\phi(H_1)=\lg u_1'\rg$, $\phi(H_2)=\lg u_2'\rg $ and $V=\perp_{i=1}^m\lg u_i',v_i'\rg  \perp \lg  w\rg $.  Then we have
$$\begin{array}{llllll}\o{H_1}&=&\{U_{1k}\di k\in\FF_q^*\},\quad   &U_{1k}&=&\lg u_1',v_1'+kw\rg\perp \lg u_2',v_2'\rg\perp (\perp_{i=3}^m\lg u_i',v_i'\rg),\\
\o{H_2}&=&\{U_{2l}\di l\in\FF_q^*\},\quad  &U_{2l}&=&\lg u_1',v_1'\rg\perp \lg u_2',v_2'+lw\rg\perp (\perp_{i=3}^m\lg u_i',v_i'\rg).\end{array}$$
Since  $$\begin{array}{ll} &U_{1k}\cap U_{2l}=\lg lv_1'+kv_2'+klw, u_1', u_2'\rg\perp (\perp_{i=3}^m\lg u_i',v_i'\rg)\\
&\rad(f\di _{U_{1k}\cap U_{2l}})=\lg ku_1'+lu_2'\rg, \, \,  Q(ku_1'+lu_2')=0,\end{array}$$
 we get  $U_{1k}\thicksim U_{2l}$,  for all $l, k$, so  that $\G(\o{H_1}\cup \o{H_2})\cong K_{14}$.

Second, assume that $\phi(H_1)$ and $\phi(H_2)$ are nonsingular. No loss, let $\phi(H_1)=\lg u_1'+v_1'\rg$ and $\phi(H_2)=\lg u_2'+v_2'\rg$.  Then we have
$$\begin{array}{llllll}\o{H_1}&=&\{U_{1k}\di k\in\FF_q^*\},\,   &U_{1k}&=&\lg u_1'+kw, v_1'+kw\rg\perp \lg u_2',v_2'\rg\perp (\perp_{i=3}^m\lg u_i',v_i'\rg),\\
\o{H_2}&=&\{U_{2l}\di l\in\FF_q^*\},\,   &U_{2l}&=&\lg u_1', v_1'\rg\perp \lg u_2'+lw, v_2'+lw\rg\perp (\perp_{i=3}^m\lg u_i',v_i'\rg),\end{array}$$
where $k=x_1^2+x_1$ and $l=x_2^2+x_2$ for some $x_1, x_2\in\FF_8^*$.
Since
$$\begin{array}{ll} &U_{1k}\cap U_{2l}=\lg lv_1'+kv_2'+klw,u_1'+v_1',u_2'+v_2'\rg\perp (\perp_{i=3}^m\lg u_i',v_i'\rg)\\
  &\rad(f\di _{U_{1k}\cap U_{2l}})=\lg k(u_1'+v_1')+l(u_2'+v_2')\rg ,\\
  &Q(k(u_1'+v_1')+l(u_2'+v_2'))=k^2+l^2\ne 0, \, \, {\rm if}\,\, k\ne l, \end{array}$$
  we know  $\G(\o{H_1}\cup\o{H_2})\cong 3K_2$.\qed

\begin{lem}\label{+3}
Suppose that $\O=\Lambda ^+$. Then   $\G=\G(\O, \D_s)$ has a H-cycle.
\end{lem}
\demo From Lemma~\ref{+1}, we just need to consider $q=8$ and $\D_s=\D_1$. That is, $W_1~W_2$ if and only if   $\dim(\rad( Q\di_{W_1\cap W_2}))=1$.
\vskip 2mm
{\bf Step 1:} {\it Find some paths for two induced subgraphs $\G(\D_1)$ and $\G(\D_2)$.}
\vskip 2mm

First  consider  $\G(\D_1)$.  Let $\G_S$ be the graph whose vertices  are all singular points of $U_0$, where  $\lg u\rg \thicksim \lg v\rg $   if and only if
$f(u,v)=0.$  Then $\G_S$ is isomorphic to  the graph, say $\G_1$ in  Lemma~\ref{esp}.(2),  by $\chi$:
$$\chi(\sum_{i=1}^m (x_iu_i+y_iv_i))=(x_2,x_1,x_3, \cdots x_m; y_2,y_1,y_3 \cdots y_n).$$
 Let $X_0, U, \b_{1i}$ and  $\b_{2i}$ have the same meaning  as in Lemma~\ref{esp}.(2).  It can be see from the proof of By Lemma~\ref{esp}.(2),   $\G_1$ contains a  H-cycle,    $\b_{2,1}\to \b_{2,k}, \b_{2,0},\b_{1,1}\to\b_{1,0},\d_{Y_1,0}$,
 where
      $$\b_{2,0}=\lg (0_m, X_0)\rg =\lg \chi(v_1)\rg ,\,
  \b_{2,1}=\lg (0_m, X_0A^\t)\rg =\lg \chi(-a_2a_0^{-1}v_2+v_3)\rg ,\,  \b_{2,k}=\chi(v_1+v_2)\rg .$$
   Thus, $\G_1$ contains a  H-path $\b_{2,1}\to \b_{2,k}$, so that
   $\G_S$ contains a H-path  of $\G_S,$
   $$P:\quad \e_0:=\lg v_1\rg, \e_1, \cdots , \e_{t-1}:=\lg v_1+v_2\rg,$$   recalling  that totally $U_0$ contains $t:=\frac{(q^m-1)(q^{m-1}+1)}{q-1}$ singular points.
  Set $\phi(H_i):=H_i\cap U_0=\lg \e_i\rg $, where $0\le i\le t-1$.
 Then by the above Lemma~\ref{+2}, $\G(\o{H_i})\cong K_7$ and   $\G(\o{H_i} \cup \o{H_{i+1}})\cong K_{14}$ where $0\le i\le t-1$.
   Set  $\FF_8^*=\{ h_1, \cdots h_{7}\} $ and $\o{H_i}=\{U_{i,h_j}\di h_j\in \FF_8^*\}$.   Then  we get a H-path of
 $\G(\D_1)$,
\begin{equation} \label{D1}  U_{0, h_1}, \cdots  U_{0, h_{7}};  \cdots  U_{t-1, h_1}, \cdots  U_{t-1, h_{7}},\end{equation}
\f where  $\bigcap_{h_j\in \FF_8^*}  U_{k, h_j}=H_k$ (see the left part of Figure 8.)

Second,  consider  $\G(\D_2)$.
Let $\G_{NS}$ be the graph whose vertices are all nonsingular points in $U_0$, where  $\lg u\rg \thicksim \lg v\rg $   if and only if
$f(u,v)=0.$ By $\chi(\G_{NS})$ we denote its image, under the isomorphism  $\chi $ in last paragraph.
First we  find a H-cycle of $\G_{NS}$.   All $\lg\rho\rg$-orbits on all nonsingular points are
$$C_Y=\{\d_{Y,k}\di 0\le k\le 8^m-2\},\, \d_{Y,k}=\lg X_0A^k,X_{1,Y}(A^k)^\t\rg,  \, \,  {\rm where }\, \,  Y\in V(m-1, 8).$$
Let $$Y_1=(1,0,0,\cdots,0), \, Y_2=(0,1,0,\cdots,0), \,  Y_0=(0,0,0,\cdots,0)\in V(m-1,8).$$ In the complete graph $\chi(\G_{NS})(\{\d_{Y,0}\di Y\})$, pick a H-path $W_0$: $\d_{Y_2,0}\to  \d_{Y_0,0}, \d_{Y_1,0}$.
Then we    get a H-cycle of $\chi(\G_{NS})$,
$$W_0,W_0^\rho,W_0^{\rho^2},\cdots,W_0^{\rho^{8^m-2}},\d_{Y_2,0}.$$
Take a H-path $\d_{Y_0,0}\to \d_{Y_1,0}$ from this cycle.     Correspondingly,  we get a
H-path  $$P':\quad \e_0'=\lg u_1+v_1\rg, \e_1', \cdots , \e_{t-1}':=\lg u_1+v_1+v_2\rg,$$ where $t'=q^{m-1}(q^m-1)$.
  Let $\phi(H_i')=\e_i'$, where $0\le i\le t'-1$.
    Then by the above Lemma~\ref{+2}, $\G(\o{H_i'})\cong \overline{K_3}$ and   $\G(\o{H_i'} \cup \o{H_{i+1}'})\cong 3K_{2}$ where $0\le i\le t'-1$.
   Set  $\FF_8=\{ h_0, h_1, \cdots h_{7}\} $ and $\o{H_i'}=\{U'_{i,h_{i,j}}\di 1\le j \le 3\}$.   Then  we get three paths  of
$\G(\D_2)$: $$U'_{0,h_{0,j}}, U'_{1, h_{1,j}}, \cdots , U'_{t'-1, h_{t'-1,j}},$$ where $1\le j \le 3$ (see the right part of Figure 8.)

\vskip 2mm
{\bf Step 2:} { \it Find a H-cycle for $\G$.}
\vskip 2mm

Recall $\phi(H_0)=\lg v_1\rg$, $\phi(H_{t-1})=\lg v_1+v_2\rg$, $\phi(H'_0)=\lg u_1+v_1\rg$ and $\phi(H_{t'-1}')=\lg u_1+v_1+v_2\rg$,  where
\vspace{-5pt}$$\begin{array}{ll}\o{H_0}=\{U_{0,h}\di h\in\FF_q^*\},\,   &U_{0,h}=\lg u_1+hw, v_1\rg\perp \lg u_2,v_2\rg\perp (\perp_{i=3}^m\lg u_i,v_i\rg),\\
\o{H_{t-1}}=\{U_{t-1,h}\di h\in\FF_q^*\},\,   &U_{t-1,h}=\lg u_1+hw, v_1+v_2\rg\perp \lg u_1+u_2, v_2\rg\perp \\
&(\perp_{i=3}^m\lg u_i,v_i\rg);\\
\o{H'_0}=\{U'_{0,k}\di k\in\FF_q^*\},\,   &U'_{0,k}=\lg u_1+kw, v_1+kw\rg\perp \lg u_2,v_2\rg\perp (\perp_{i=3}^m\lg u_i,v_i\rg),\\
\o{H'_{t'-1}}=\{U_{t'-1,l}\di l\in\FF_q^*\},\,   &U'_{t'-1,l}=\lg u_1+lw, v_1+v_2+lw\rg\perp \lg u_1+u_2, v_2\rg\perp \\
&(\perp_{i=3}^m\lg u_i,v_i\rg),\end{array}$$
where  $k=x_1^2+x_1$ and $l=x_2^2+x_2$  for some $x_1, x_2\in\FF_8^*$  so that
 $k,l\in \{\th,\th^2,\th^4\}$, provided $\FF_8=\FF_2(\th )$  where   $\th^3+\th+1=0.$

Now we are finding  the condition of $U_{0,h}\thicksim U'_{0,k}$.
Since  $$U_{0,h}\cap U'_{0,k}=\lg u_1+av_1+hw\rg\perp(\perp_{i=2}^m\lg u_i,v_i\rg),$$
for some $a\in \FF_8$, we get $\rad(f\di_{U_{0,h}\cap U'_{0,k}})=\lg u_1+av_1+hw\rg$.  So  $U_{0,h}\thicksim U'_{0,k}$ if and only if
$Q(u_1+av_1+hw)=0$, that is $h^2+a=0$. Noting that  $\lg u_1+av_1+hw\rg=\lg u_1+av_1+k(1+a)w\rg$, we have
$h=k(1+a)$, where  $h, a \in\FF_8^*$ and $k\in \{ \th,\th^2,\th^4\}.$ Therefore,
 all  solutions for  $(a,h,k)$ are
 $$(\th,\th^4,\th), (\th^6,\th^3,\th), (\th^2,\th,\th^2), (\th^5,\th^6,\th^2), (\th^3,\th^5,\th^4), (\th^4,\th^2,\th^4).$$
 Thus by picking three pairs  $h_1=\th^4,k_1=\th$; $h_2=\th,k_2=\th^2$;  and $h_3=\th^5,k_3=\th^4$,  we have $U_{0,h_i}\thicksim U'_{0,k_i}$, where $1\le i \le 3$.

Recall $\phi(H_{t-1})=\lg v_1+v_2\rg$ and $\phi(H'_{t-1})=\lg u_1+v_1+v_2\rg $. Similar to last paragraph, we have three pairs $(h_i, l_i)$ such that  $U_{t-1,h_i}\thicksim U'_{t'-1,l_i}$, where $1\le i \le 3$.

Simply,  write the above three H-paths of $\G(\D_2)$ as  $U'_{0,k_i}\to U'_{t'-1,l_i}$, where $1\le i\le 3$.
 By $P_1$ we denote the H-path of   $\G(\D_1)\setminus \{U_{0,h_i},U_{t-1,h_i}\di 1\le i \le 3\})$ derived from    the H-path in Eq(\ref{D1}) of  $\G(\D_1)$ by removing these six points.
 Noting  $\G(\o{H_0}\cup\o{H_{t-1}})\cong K_{14}$,  we finally get a H-cycle $C$ for $\G$ (see Figure 8.):
$$\begin{array}{ll} C:\,& U_0,U_{0,h_1},U'_{0,k_1}\to U'_{t'-1,l_1},U_{t-1,h_1},U_{0,h_2},U'_{0,k_2}\to U'_{t'-1,l_2},U_{t-1,h_2},U_{0,h_3},\\
& U'_{0,k_3}\to U'_{t'-1,l_3},U_{t-1,h_3}, P_1, U_0.\end{array}$$

\begin{figure}[htbp]
  \centering
  \begin{tikzpicture}[
    scale=0.75, transform shape,
    p/.style={draw, circle, fill=white, inner sep=0pt, minimum size=2.6mm},
    d/.style={fill=white, inner sep=2pt},
    box/.style={draw, rounded corners=6pt},
    lab/.style={anchor=north, yshift=3pt, font=\footnotesize},
    tlab/.style={anchor=south, font=\small},
    midarrow/.style={
      postaction=decorate,
      decoration={markings, mark=at position 0.5 with {\arrow{angle 45}}}
    },
    remidarrow/.style={
      postaction=decorate,
      decoration={markings, mark=at position 0.5 with {\arrowreversed{angle 45}}}
    }
  ]
    \node[p] (a10)  at (0, 0)  {};
    \node[d] (d10)  at (1, 0)  {$\cdots$};
    \node[p] (a20)  at (2, 0)  {};
    \node[p] (a30)  at (4, 0)  {};
    \node[d] (d20)  at (5, 0)  {$\cdots$};
    \node[p] (a40)  at (6, 0)  {};
    \node[p] (a50)  at (7, 0)  {};
    \node[p] (a60)  at (8, 0)  {};
    \node[p] (a70)  at (9, 0)  {};

    \node[p] (b10)  at (12, 0)  {};
    \node[p] (b20)  at (13, 0)  {};
    \node[p] (b30)  at (14, 0)  {};
    \draw[remidarrow] (a10.east) -- (d10.west);
    \draw[remidarrow] (d10.east) -- (a20.west);
    \draw[remidarrow] (a20.east) -- (a30.west);
    \draw[remidarrow] (a30.east) -- (d20.west);
    \draw[remidarrow] (d20.east) -- (a40.west);
    \draw[midarrow] (a70.east) -- (b10.west);

    \node[d] (v1)  at (1, -1)  {$\vdots$};
    \node[d] (v2)  at (12, -1)  {$\vdots$};
    \node[d] (v3)  at (13, -1)  {$\vdots$};
    \node[d] (v4)  at (14, -1)  {$\vdots$};
    \node[p] (a)  at (10.5, -1)  {};

    \node[p] (a11)  at (0, -2)  {};
    \node[d] (d11)  at (1, -2)  {$\cdots$};
    \node[p] (a21)  at (2, -2)  {};
    \node[p] (a31)  at (4, -2)  {};
    \node[d] (d21)  at (5, -2)  {$\cdots$};
    \node[p] (a41)  at (6, -2)  {};
    \node[p] (a51)  at (7, -2)  {};
    \node[p] (a61)  at (8, -2)  {};
    \node[p] (a71)  at (9, -2)  {};

    \node[p] (b11)  at (12, -2)  {};
    \node[p] (b21)  at (13, -2)  {};
    \node[p] (b31)  at (14, -2)  {};
    \draw[remidarrow] (a11.east) -- (d11.west);
    \draw[remidarrow] (d11.east) -- (a21.west);

    \draw[midarrow] (a31.east) -- (d21.west);
    \draw[midarrow] (d21.east) -- (a41.west);
    \draw[remidarrow] (a71.east) -- (b11.west);

    \draw[midarrow] (a.north) -- (a70.south);
    \draw[midarrow] (b10.south) -- (v2.north);
    \draw[midarrow] (v2.south) -- (b11.north);
    \draw[midarrow] (a71.north) -- (a60.south);
    \draw[midarrow, rounded corners=6pt]
    (a60.north) -- (8,0.9) -- (13,0.9) -- (b20.north);
    \draw[midarrow] (b20.south) -- (v3.north);
    \draw[midarrow] (v3.south) -- (b21.north);
    \draw[midarrow, rounded corners=6pt]
    (b21.south) -- (13,-2.9) -- (8,-2.9) -- (a61.south);
    \draw[midarrow] (a61.north) -- (a50.south);

    \draw[midarrow, rounded corners=6pt]
    (a50.north) -- (7,1.1) -- (14,1.1) -- (b30.north);
    \draw[midarrow] (b30.south) -- (v4.north);
    \draw[midarrow] (v4.south) -- (b31.north);
    \draw[midarrow, rounded corners=6pt]
    (b31.south) -- (14,-3.1) -- (7,-3.1) -- (a51.south);
    \draw[midarrow] (a51.north) -- (a40.south);

    \draw[midarrow] (a10.south) -- ($(a10)!0.36!(a21)$);
    \draw[midarrow] ($(a10)!0.64!(a21)$) -- (a21.north);

    \draw[midarrow, rounded corners=6pt]
    (a11.south) -- (0,-2.9) -- (4,-2.9) -- (a31.south);
    \draw[midarrow] (a41.north) -- (a.south);

    \draw[box] (-0.7,0.7) rectangle (9.7,-2.7);
    \draw[box] (-0.5,0.5) rectangle (2.5,-0.55);
    \draw[box] (-0.5,-1.5) rectangle (2.5,-2.55);
    \draw[box] (3.5,0.5) rectangle (9.5,-0.55);
    \draw[box] (3.5,-1.5) rectangle (9.5,-2.55);

    \draw[box] (11.3,0.7) rectangle (14.7,-2.7);
    \draw[box] (11.5,0.5) rectangle (14.5,-0.55);
    \draw[box] (11.5,-1.5) rectangle (14.5,-2.55);

    \node[lab] at (a.south) {$U_{0}$};
    \node[lab] at (a10.south) {$U_{1,h_7}$};
    \node[lab] at (a20.south) {$U_{1,h_1}$};
    \node[lab] at (a30.south) {$U_{0,h_7}$};
    \node[lab] at (a40.south) {$U_{0,h_4}$};
    \node[lab] at (a50.south) {$U_{0,h_3}$};
    \node[lab] at (a60.south) {$U_{0,h_2}$};
    \node[lab] at (a70.south) {$U_{0,h_1}$};
    \node[lab] at (b10.south) {$U'_{0,k_1}$};
    \node[lab] at (b20.south) {$U'_{0,k_2}$};
    \node[lab] at (b30.south) {$U'_{0,k_3}$};

    \node[lab] at (a11.south) {$U_{t-2,h_7}$};
    \node[lab] at (a21.south) {$U_{t-2,h_1}$};
    \node[lab] at (a31.south) {$U_{t-1,h_4}$};
    \node[lab] at (a41.south) {$U_{t-1,h_7}$};
    \node[lab] at (a51.south) {$U_{t-1,h_3}$};
    \node[lab] at (a61.south) {$U_{t-1,h_2}$};
    \node[lab] at (a71.south) {$U_{t-1,h_1}$};
    \node[lab] at (b11.south) {$U'_{t'-1,l_1}$};
    \node[lab] at (b21.south) {$U'_{t'-1,l_2}$};
    \node[lab] at (b31.south) {$U'_{t'-1,l_3}$};

    \node[tlab] at (0,0.6) {$\G(S)$};
    \node[tlab] at (15,0.6) {$\G(NS)$};

    \node[tlab] at (2,-1.15) {$\o{H_1}$};
    \node[tlab] at (4,-1.15) {$\o{H_0}$};

    \node[tlab] at (2,-1.65) {$\o{H_{t-2}}$};
    \node[tlab] at (4,-1.65) {$\o{H_{t-1}}$};

    \node[tlab] at (12.5,-1.15) {$\o{H'_0}$};
    \node[tlab] at (12.5,-1.65) {$\o{H'_{t'-1}}$};

  \end{tikzpicture}
  \caption{{\footnotesize H-cyce $C$ of the graph $\G$ in Lemma 7.2}}
\end{figure}
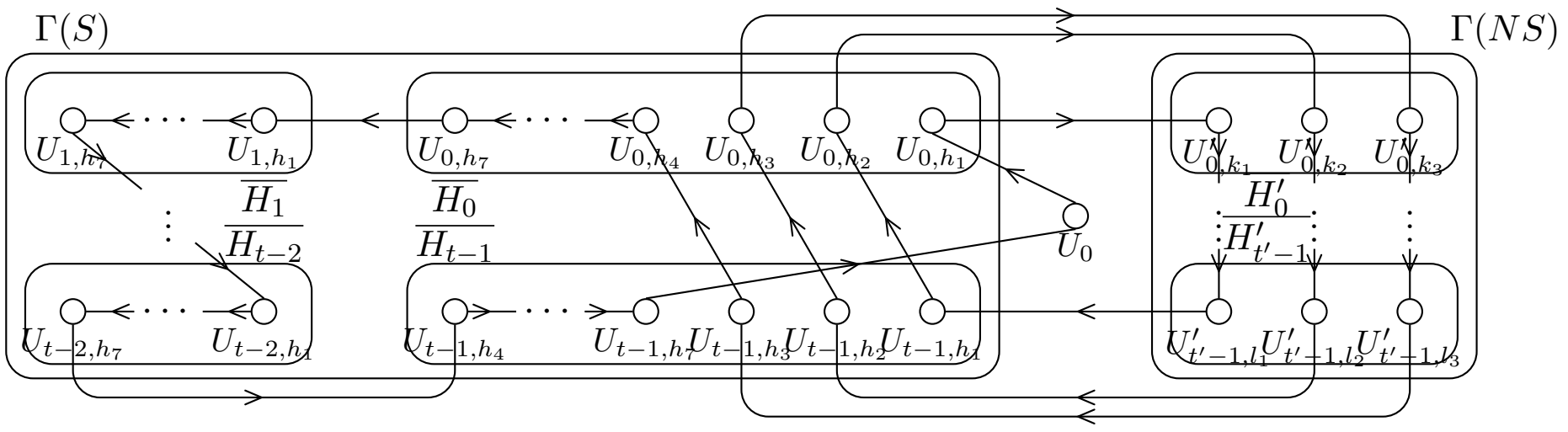

\subsection{ $\O=\Lambda ^-$}
Let    $U_1=\perp_{i=1}^{m-1}\lg u_i,v_i\rg\perp\lg u_m+\ld_m w,v_m+\mu_m w\rg \in \O:=\Lambda ^-$.
With the same arguments as in Lemmas~\ref{+1} and \ref{+2}, respectively, we have the following two lemmas.
\begin{lem}\label{-1} The stabilizer $M_{U_1}$ has two nontrivial orbits on $\O$:
$$\begin{array}{lll}
&\D_1=\{W\in\O \di \dim(\rad (Q\di_{U_1\cap W})=0)\},\quad &|\D_1|=(q^m-1)(q^{m-1}+1);\\

&\D_2=\{W\in\O \di \dim(\rad (Q\di_{U_1\cap W})=1)\}, \quad  &|\D_2|=(\frac q2-1)\cdot q^{m-1}(q^m-1).\end{array}$$
\f If $q=4$, then $|\D_1|, |\D_2|\ge \frac 13|\O|$. \end{lem}
\begin{lem}\label{-2}
Let $H_1$ and $H_2$ be  two maximal subspaces of $U_0$ such that $\phi(H_1)$ and $\phi(H_2)$ are two  singular points (resp. nonsingular points) which are mutually orthogonal. Then $\G(\o{H_1} \cup \o{H_2})$ is $K_{14}$ (resp.  $3K_2$).
\end{lem}

\begin{lem}\label{-3}
Suppose that $\O=\Lambda ^-$. Then   $\G=\G(\O, \D_s)$ has a H-cycle.
\end{lem}
\demo For our $U_1$, let $\G_S$ and $\G_{NS}$ have the same meaning in the proof of Lemma~\ref{+3}.

\vskip 2mm
{\bf Step 1:} {\it Find some paths for two induced subgraphs $\G(\D_1)$ and $\G(\D_2)$.}
\vskip 2mm
Similar to Lemma~\ref{+3},  one may prove that  $\G_S$ contains a H-path.

To find a H-path of $\G_{NS}$, let $A,\rho,r,0_{m-1},X_0$ and $X_{r,Y}$ have the same meanings in Lemma~$\ref{esm}.(2)$. Set
$$\begin{array}{llll}
\{\b\},& &  &\b=\lg(0_{m-1},0_{m-1},0,1)\rg .\\
\{\b_{c}\}  &&  &\b_{c}=\lg(0_{m-1},0_{m-1},1,c)\rg ,\,  c\in\FF_8;\\
C_{Y}&=&\{\d_{Y,k} \di 0\le k \le 8^{m-1}-2\},\quad &\d_{Y,k}=\lg(X_0A^k,X_{1,Y}(A^k)^\t,0,0)\rg,\\
&&& Y\in V(m-2,8);\\
D&=&\{\g_k \di 0\le k \le 8^{m-1}-2\},\quad &\g_k=\lg(0_{m-1},X_0(A^k)^\t,0,1)\rg,\\
D_{b,Y}&=&\{\g_{b,Y;k} \di 0\le k \le 8^{m-1}-2\},\quad &\g_{b,Y;k}=\lg(X_0A^k,X_{b,Y}(A^k)^\t,0,1)\rg,\\
&&& Y\in V(m-2,8), b\ne 1\in\FF_8;\\
E_c&=&\{\ph_{c,k} \di 0 \le k \le 8^{m-1}-2\},\quad &\ph_{c,k}=\lg(0_{m-1},X_0(A^k)^\t,1,c)\rg\\
&&& c\in \FF_8;\\
E_{c,d,Y}&=&\{\ph_{c,d,Y;k} \di 0 \le k \le 8^{m-1}-2\},\quad &\ph_{c,d,Y;k}=\lg(X_0A^k,X_{d,Y}(A^k)^\t,1,c)\rg\\
&&& Y\in V(m-2,8), c\in \FF_8,\\
&&&d\in \FF_8\setminus \{1+c+c^2\}.
\end{array}$$
Then  $\O$ is the union of disjoints of the all $\lg\rho\rg$-orbits, that is
$$\O=\{\b\}\cup \{\b_{c}\}\cup (\cup _{Y}C_Y)\cup D\cup (\cup _{b, Y}D_{b, Y}) \cup (\cup_c E_c)\cup (\cup _{c,d, Y}E_{c,d,Y}).\vspace{-2mm}$$
Moreover,   $\G_N(\{\d_{Y,k}\di Y\})$ for given $k$, $\G_N(\{\g_{b,Y;k}\di Y\})$ for given $b,k$, and $\G_N(\{\ph_{c,d,Y;k}\di Y\})$ for given $c,d,k$ are complete graphs.
Note that $\b$ is adjacent all points in $C_Y$, $D_{b,Y}$ and $D$;  and $\b_c$ is adjacent to all points in $C_Y$, $E_{c,d,Y}$ and $E_c$.

Let $$\FF_{8}=\{b_i\di 0\le i \le 7\}=\{c_i\di 0\le i \le 7\}=\{d_i\di 0\le i \le 7\}$$
 where $b_0=1, c_0=0, c_1=1$ and $d_0=1+c+c^2$ for given $c$ in $\ph_{c,d,Y;k}$.
Since $\ph_{c,k}$ is adjacent to all points in $\{\g_{b,Y;k}\di Y\}$ for given $k$ and $\g_k$ is adjacent to all points in $\{\ph_{c,d,Y;k}\di Y\}$ for given $k$, we can let $\ph_{c_i,k}\thicksim \g_{b_i,Y;k}$, where $1\le i \le q-1$, $\g_k\thicksim \ph_{c_1=1,d_1,Y;k}$ and $\ph_{c_0=0,k}\thicksim \ph_{c_1=1,d_2,Y;k}$.

If  $m=3$, we can find a H-cycle for $\G$ by Magma, directly.
Suppose $m\ge 4$.  Take ectors in  $V(m-2,q)$
$$\begin{array}{ll} &Y_0=0_{m-2},\quad   Y_1=(1,0,0,\cdots,0), \quad   Y_2=(0,1,0,\cdots,0),  \\
  & Y_3(c_0)=0_{m-2}, \quad Y_3(c_i)=(c_{i-1}+c_i,0,0,\cdots,0),  (1\le i \le q-1), \end{array} $$

Pick a H-path $W_0$: $\d_{Y_2,0},\cdots,\d_{Y_0,0},\d_{Y_1,0}$ in the complete graph $\G_N(\{\d_{Y,0}\di Y\})$. Since   $\d_{Y_1,0}\thicksim \d_{Y_2,1}$,  we get a H-cycle for $\G_N(\{C_Y\di Y\})$:
$$W=W_0,W_0^\rho,W_0^{\rho^2},\cdots,W_0^{\rho^{8^{m-1}-2}},\d_{Y_2,0}.$$

Pick a H-path $U_{0}(b)$: $\g_{b,Y_2;0},\cdots,\g_{b,Y_0,0},\g_{b,Y_1;0}$ in the complete graph $\G_N(\{\g_{b,Y;0}\di Y\}\cup\{\ph_{c,0})\}$. Since  $\g_{b,Y_1;0}\thicksim\g_{b,Y_2;1}$,  we get a H-cycle for $\G_N(\{D_{b,Y}\di b\ne b_0,Y\}\cup\{E_c\di c\ne c_0\}):$
$$U=U_1(b_1),U_1(b_1)^\rho,\cdots,U_1(b_1)^{\rho^{8^{m-1}-2}},U_1(b_2),\cdots,U_1(b_{q-1})^{\rho^{8^{m-1}-2}},\ph_{b_{1},Y_2;0}.$$

Pick a H-path $R_0(d)$: $\ph_{c_1,d,Y_2;0},\cdots,\ph_{c_1,d,Y_3(c_1);0},\ph_{c_1,d,Y_1;0}$
 in the complete graphs $\G_N(\{\ph_{c_1,d_1,Y;0}\di Y\}\cup\{\g_0\})$, if $d=d_1$;  $\G_N(\{\ph_{c_1,d_2,Y;0}\di Y\}\cup\{\ph_{c_0,0}\})$ if $d_2$ and $\G_N(\{\ph_{c_1,d,Y;0}\di Y\})$ if $d\in \FF_8\setminus
 \{d_0,  d_1, d_2\}$. Since $\ph_{c,d,Y_1;0}\thicksim \ph_{c,d,Y_2;1}$,  we get a H-cycle for $\G_N(\{E_{c_1,d,Y;0}\di d\ne d_0 , Y\}\cup\{D,E_{c_0}\})$,
$$R(c_1)=R_0(d_1),R_0(d_1)^\rho,\cdots,R_0(d_1)^{\rho^{8^{m-1}-2}},R_0(d_2),\cdots,R_0(d_{q-1})^{\rho^{8^{m-1}-2}},\ph_{c_1,d_1,Y_2;0}.$$
Similarly, we can get all H-cycle $R(c)$ for $\G(\{E_{c,d,Y}\di d\ne 1+c+c^2,Y\})$, where $c\ne c_1$.

(i) Inserting $\{\b_{c_i}\di i\ne q-2\}$ to  the H-cycle $W$ and get a H-path $\d_{Y_1,0},\b_{c_0},\cdots,\b_{c_{q-2}}\to\d_{Y_0,0}$;  (ii) take a
  a H-path $\g_{b_1,Y_1;0}\to\g_{b_1,Y_0;0}$ from the H-cycle  $U$; (iii) take  a H-path $\ph_{c_0,d_1,Y_1;1}\to\ph_{c_0,d_1,Y_0;1}$  from the H-cycle  $R(c_0)$;  and a H-path $\ph_{c_i,d_1,Y_3(c_i);i+1}\to\ph_{c_i,d_1,Y_0;i+1}$ from the H-cycle $R(c_i)$, where $1\le i \le q-1$.
Noting that $\ph_{c_i,d_1,Y_0;i+1}\thicksim\ph_{c_{i+1},d_1,Y_3(c_{i+1});i+2}$,  we finally get  a H-cycle for $\G_N$:
$$\begin{array}{l}\d_{Y_1,0},\b_{c_0},\cdots,\b_{c_{q-2}}\to\d_{Y_0,0},\b,\g_{b_1,Y_1;0}\to\g_{b_1,Y_0;0},\ph_{c_0,d_1,Y_1;1}\to\ph_{c_0,d_1,Y_0;1},\\
\ph_{c_1,d_1,Y_3(c_1);2}\to\ph_{c_1,d_1,Y_0;2},\ph_{c_2,d_1,Y_3(c_2);3},\cdots,\ph_{c_{q-1},d_1,Y_0;q},\b_{c_{q-1}},\d_{Y_1,0},
\end{array}$$
as desired.

\vskip 2mm
{\bf Step 2:} { \it Find a H-cycle for $\G$.}
\vskip 2mm
Completely same with the Step 2 of Lemma ~\ref{+3}
\qed

\end{document}